\documentclass[letterpaper,11pt]{article}
\usepackage{amsmath, amsthm, amssymb}    % May not all be necessary
\usepackage{bridges}			% Custom bridges proceedings style
\usepackage{graphicx}			% For including pictures
\usepackage[colorlinks=true, urlcolor=blue, citecolor=black, linkcolor=black]{hyperref}  
\usepackage{subcaption}			% For captioning multi-panel figures

\usepackage{graphicx}
\usepackage{tikz}
\usetikzlibrary{arrows}
\usetikzlibrary{decorations.markings}
\usetikzlibrary{decorations.pathreplacing}
\usetikzlibrary{patterns}
\usetikzlibrary{shapes.geometric}
\usetikzlibrary{matrix}
\usepackage{tikz-3dplot}
\usepackage{tkz-graph}
\usepackage{tikz-cd}

\usepackage{xcolor}
\usepackage{color}
\usepackage{visualalgebra}  %% Put this *after* the TikZ packages

\title{Revealing the Hidden Beauty of Finite Groups with Cayley Graphs}

\author{Matthew Macauley\textsuperscript{1} \vspace{10pt}\\
\textsuperscript{1}School of Mathematical \& Statistical Sciences, Clemson University; macaule@clemson.edu} % end \author
\date{}					% Suppress any date on submissions

\begin{document}

\maketitle

% Prevent page number 1 from being printed on the first page.
\thispagestyle{empty}

\begin{abstract}
Group theory involves the study of symmetry, and its inherent beauty gives it the potential to be one of the most accessible and enjoyable areas of mathematics, for students and non-mathematicians alike. Unfortunately, many students never get a glimpse into the more alluring parts of this field because ``traditional'' algebra classes are often taught in a dry axiomatic fashion, devoid of visuals. This article will showcase aesthetic pictures that can bring this subject to life. It will also leave the reader with some (intentionally) unanswered puzzles that undergraduate students, hobbyists, and mathematical artists can explore and answer, and even create new versions themselves. 
\end{abstract}

%%%%%%%%%%%%%%%%%%%%%%%%%%%%%%%%%%%%%%%%%%%
\section*{Cayley Graphs and Tables}

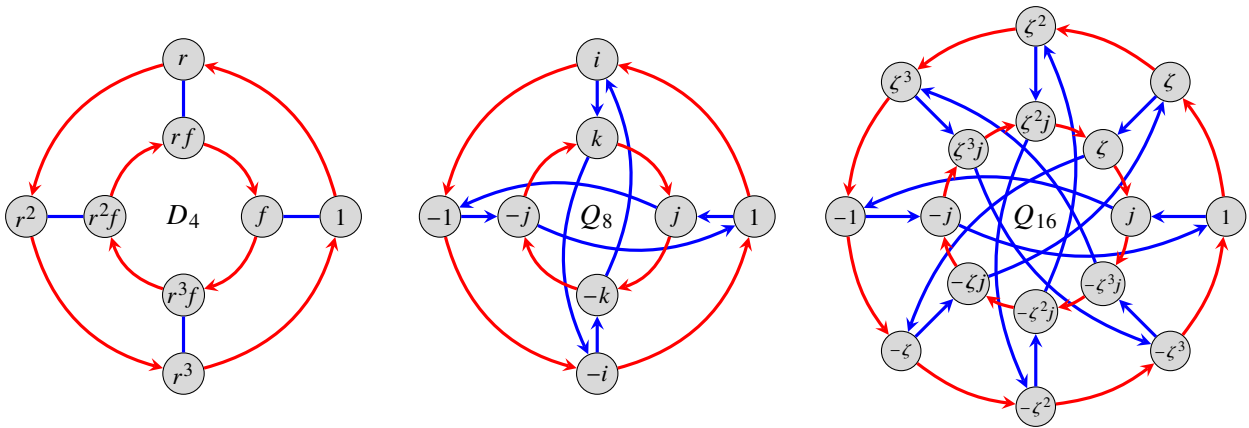
\begin{figure}[!ht]
  \tikzstyle{v} = [circle, draw, fill=lightgray,inner sep=0pt,
    minimum size=5.5mm]
\begin{tikzpicture}[scale=1.1]
  \tikzstyle{R-out} = [draw, very thick, eRed,-stealth,bend right=15]
  \tikzstyle{R-in} = [draw, very thick, eRed,-stealth,bend left=12]
  \tikzstyle{B} = [draw, very thick, eBlue,-stealth,bend right=25]
  \begin{scope}[shift={(0,0)},scale=.95]
      \tikzstyle{every node}=[font=\footnotesize]
      \node (1) at (0:2) [v] {$1$};
      \node (r) at (90:2) [v] {$r$};
      \node (r2) at (180:2) [v] {$r^2$};
      \node (r3) at (270:2) [v] {$r^3$};
      \node (s) at (0:1) [v] {$f$};
      \node (r3s) at (270:1) [v] {$r^3\!f$};
      \node (r2s) at (180:1) [v] {$r^2\!f$};
      \node (rs) at (90:1) [v] {$rf$};
      \draw [bb] (1) to (s); 
      \draw [bb] (r2) to (r2s); 
      \draw [bb] (r) to (rs); 
      \draw [bb] (r3) to (r3s); 
      \draw [r] (1) to [bend right] (r);
      \draw [r] (r) to [bend right] (r2);
      \draw [r] (r2) to [bend right] (r3);
      \draw [r] (r3) to [bend right] (1);
      \draw [r] (s) to [bend left=28] (r3s);
      \draw [r] (r3s) to [bend left=28] (r2s);
      \draw [r] (r2s) to [bend left=28] (rs);
      \draw [r] (rs) to [bend left=28] (s);
      \node at (0,0) {\normalsize $D_4$};
    \end{scope}
    \begin{scope}[shift={(5,0)},scale=.95]
      \tikzstyle{every node}=[font=\footnotesize]
      \node (1) at (0:2) [v] {$1$};
      \node (r) at (90:2) [v] {$i$};
      \node (r2) at (180:2) [v] {$-1$};
      \node (r3) at (270:2) [v] {$-i$};
      \node (s) at (0:1) [v] {$j$};
      \node (r3s) at (270:1) [v] {$-k$};
      \node (r2s) at (180:1) [v] {$-j$};
      \node (rs) at (90:1) [v] {$k$};
      \draw [b] (1) to (s); \draw [B] (s) to (r2);
      \draw [b] (r2) to (r2s); \draw [B] (r2s) to (1);
      \draw [b] (r) to (rs); \draw [B] (rs) to (r3);
      \draw [b] (r3) to (r3s); \draw [B] (r3s) to (r);
      \draw [r] (1) to [bend right] (r);
      \draw [r] (r) to [bend right] (r2);
      \draw [r] (r2) to [bend right] (r3);
      \draw [r] (r3) to [bend right] (1);
      \draw [r] (s) to [bend left=28] (r3s);
      \draw [r] (r3s) to [bend left=28] (r2s);
      \draw [r] (r2s) to [bend left=28] (rs);
      \draw [r] (rs) to [bend left=28] (s);
      \node at (0,0) {\normalsize $Q_8$};
    \end{scope}
    \begin{scope}[shift={(10.3,0)},scale=1.15,auto]
    \tikzstyle{v} = [circle, draw, fill=lightgray,inner sep=0pt,
    minimum size=5.25mm]
    \tikzstyle{every node}=[font=\scriptsize]
      \node (s) at (0:1) [v] {$j$};
      \node (rs) at (45:1) [v] {$\zeta$};
      \node (r2s) at (90:1) [v] {$\zeta^2\!j$};
      \node (r3s) at (135:1) [v] {$\zeta^3\!j$};
      \node (r4s) at (180:1) [v] {$-j$};
      \node (r5s) at (225:1) [v] {$-\zeta\!j$};
      \node (r6s) at (270:1) [v] {\tiny $-\zeta^2\!j$};
      \node (r7s) at (315:1) [v] {\tiny$-\zeta^3\!j$};
      \node (1) at (0:2) [v] {$1$};
      \node (r) at (45:2) [v] {$\zeta$};
      \node (r2) at (90:2) [v] {$\zeta^2$};
      \node (r3) at (135:2) [v] {$\zeta^3$};
      \node (r4) at (180:2) [v] {$-1$};
      \node (r5) at (225:2) [v] {\tiny $-\zeta$};
      \node (r6) at (270:2) [v] {\tiny $-\zeta^2$};
      \node (r7) at (315:2) [v] {\tiny $-\zeta^3$};
      \draw [b] (1) to (s); \draw [B] (s) to (r4);
      \draw [b] (r4) to (r4s); \draw [B] (r4s) to (1);
      \draw [b] (r) to (rs); \draw [B] (rs) to (r5);
      \draw [b] (r5) to (r5s); \draw [B] (r5s) to (r);
      \draw [b] (r2) to (r2s); \draw [B] (r2s) to (r6);
      \draw [b] (r6) to (r6s); \draw [B] (r6s) to (r2);
      \draw [b] (r3) to (r3s); \draw [B] (r3s) to (r7);
      \draw [b] (r7) to (r7s); \draw [B] (r7s) to (r3);
      \draw [R-out] (1) to (r);
      \draw [R-out] (r) to (r2);
      \draw [R-out] (r2) to (r3);
      \draw [R-out] (r3) to (r4);
      \draw [R-out] (r4) to (r5);
      \draw [R-out] (r5) to (r6);
      \draw [R-out] (r6) to (r7);
      \draw [R-out] (r7) to (1);
      \draw [R-in] (s) to (r7s);
      \draw [R-in] (r7s) to (r6s);
      \draw [R-in] (r6s) to (r5s);
      \draw [R-in] (r5s) to (r4s);
      \draw [R-in] (r4s) to (r3s);
      \draw [R-in] (r3s) to (r2s);
      \draw [R-in] (r2s) to (rs);
      \draw [R-in] (rs) to (s);
      \node at (0,0) {\normalsize $Q_{16}$};
    \end{scope}   
    \end{tikzpicture}
    \caption{Several Cayley graphs. Multiplication of elements will be read from left-to-right.}\label{fig:D4-Q8-Q16}
\end{figure}

Every undergraduate algebra student learns about the dihedral and quaternion groups, which are typically presented as
\[
D_n=\big\<r,f\mid r^n=f^2=1,\,rfr=f\big\>,\qquad
Q_8=\big\<i,j,k\mid i^2=j^2=k^2=ijk=-1\big\>.
\]
These presentations can be visualized by a \emph{Cayley graph}, where the nodes represent elements, and the arrows represent generators. Examples of these for $D_4$ and $Q_8$ are shown in Figure~\ref{fig:D4-Q8-Q16}, alongside a third, which describes a \emph{generalized quaternion group}. This group, denoted $Q_{16}$, is defined by replacing the fourth root of unity $i=e^{2\pi i/4}$ in $Q_8$ with an eighth root of unity, $\zeta=e^{2\pi i/8}=\frac{\sqrt{2}}{2}+\frac{\sqrt{2}}{2}i$.

    %% Generalized quaternion group (dicyclic) of order 12 and 16.

In Figure~\ref{fig:D4-Q8-Q16} and henceforth, edges that are undirected are assumed to be bidirected; these represent generators of order $2$. Additionally, multiplication will be read from left-to-right. For example, in $D_4$, the red arrow represents $r$, and the blue arrow, $f$. Thus, $rf$ represents the red-blue path. Specifically, for any $g\in D_4$, the element $g\cdot (rf)$ can be found by following the red-blue path \emph{from} node $g$.

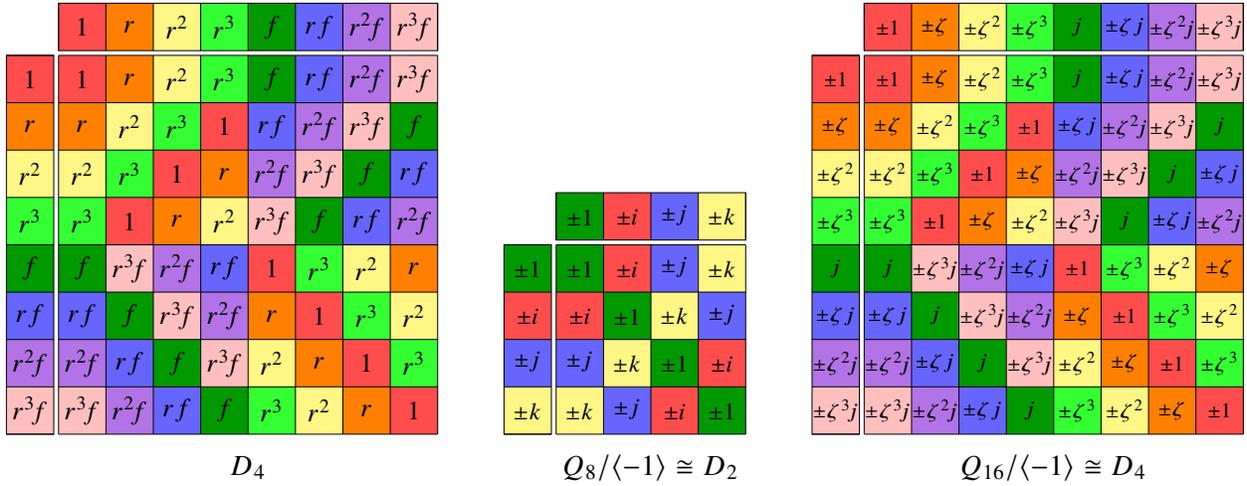
\begin{figure}[!ht]
  \begin{tikzpicture}[scale=0.63]
    \colorlet{color-1}{tRed}
    \colorlet{color-r}{tOrange}
    \colorlet{color-r2}{tYellow}
    \colorlet{color-r3}{tLime}
    \colorlet{color-f}{tGreen!80!black}
    \colorlet{color-rf}{tBlue}
    \colorlet{color-r2f}{tPurple}
    \colorlet{color-r3f}{tPink}
    \newcommand*{\n}{9}%
    \begin{scope}[shift={(0,0)}]
    \tikzstyle{every node}=[font=\small]
    \path[fill=color-1] (-.1,8) rectangle ++(1,1);
    \path[fill=color-r] (-.1,7) rectangle ++(1,1);
    \path[fill=color-r2] (-.1,6) rectangle ++(1,1);
    \path[fill=color-r3] (-.1,5) rectangle ++(1,1);
    \path[fill=color-f] (-.1,4) rectangle ++(1,1);
    \path[fill=color-rf] (-.1,3) rectangle ++(1,1);
    \path[fill=color-r2f] (-.1,2) rectangle ++(1,1);
    \path[fill=color-r3f] (-.1,1) rectangle ++(1,1);
    \path[fill=color-1] (1,9.1) rectangle ++(1,1);
    \path[fill=color-r] (2,9.1) rectangle ++(1,1);
    \path[fill=color-r2] (3,9.1) rectangle ++(1,1);
    \path[fill=color-r3] (4,9.1) rectangle ++(1,1);
    \path[fill=color-f] (5,9.1) rectangle ++(1,1);
    \path[fill=color-rf] (6,9.1) rectangle ++(1,1);
    \path[fill=color-r2f] (7,9.1) rectangle ++(1,1);
    \path[fill=color-r3f] (8,9.1) rectangle ++(1,1);
    \path[fill=color-1] (1,8) rectangle ++(1,1);
    \path[fill=color-1] (2,5) rectangle ++(1,1);
    \path[fill=color-1] (3,6) rectangle ++(1,1);
    \path[fill=color-1] (4,7) rectangle ++(1,1);
    \path[fill=color-1] (5,4) rectangle ++(1,1);
    \path[fill=color-1] (6,3) rectangle ++(1,1);
    \path[fill=color-1] (7,2) rectangle ++(1,1);
    \path[fill=color-1] (8,1) rectangle ++(1,1);
    \path[fill=color-r] (1,7) rectangle ++(1,1);
    \path[fill=color-r] (2,8) rectangle ++(1,1);
    \path[fill=color-r] (3,5) rectangle ++(1,1);
    \path[fill=color-r] (4,6) rectangle ++(1,1);
    \path[fill=color-r] (5,3) rectangle ++(1,1);
    \path[fill=color-r] (6,2) rectangle ++(1,1);
    \path[fill=color-r] (7,1) rectangle ++(1,1);
    \path[fill=color-r] (8,4) rectangle ++(1,1);
    \path[fill=color-r2] (1,6) rectangle ++(1,1);
    \path[fill=color-r2] (2,7) rectangle ++(1,1);
    \path[fill=color-r2] (3,8) rectangle ++(1,1);
    \path[fill=color-r2] (4,5) rectangle ++(1,1);
    \path[fill=color-r2] (5,2) rectangle ++(1,1);
    \path[fill=color-r2] (6,1) rectangle ++(1,1);
    \path[fill=color-r2] (7,4) rectangle ++(1,1);
    \path[fill=color-r2] (8,3) rectangle ++(1,1);
    \path[fill=color-r3] (1,5) rectangle ++(1,1);
    \path[fill=color-r3] (2,6) rectangle ++(1,1);
    \path[fill=color-r3] (3,7) rectangle ++(1,1);
    \path[fill=color-r3] (4,8) rectangle ++(1,1);
    \path[fill=color-r3] (5,1) rectangle ++(1,1);
    \path[fill=color-r3] (6,4) rectangle ++(1,1);
    \path[fill=color-r3] (7,3) rectangle ++(1,1);
    \path[fill=color-r3] (8,2) rectangle ++(1,1);
    \path[fill=color-f] (1,4) rectangle ++(1,1);
    \path[fill=color-f] (2,3) rectangle ++(1,1);
    \path[fill=color-f] (3,2) rectangle ++(1,1);
    \path[fill=color-f] (4,1) rectangle ++(1,1);
    \path[fill=color-f] (5,8) rectangle ++(1,1);
    \path[fill=color-f] (6,5) rectangle ++(1,1);
    \path[fill=color-f] (7,6) rectangle ++(1,1);
    \path[fill=color-f] (8,7) rectangle ++(1,1);
    \path[fill=color-rf] (1,3) rectangle ++(1,1);
    \path[fill=color-rf] (2,2) rectangle ++(1,1);
    \path[fill=color-rf] (3,1) rectangle ++(1,1);
    \path[fill=color-rf] (4,4) rectangle ++(1,1);
    \path[fill=color-rf] (5,7) rectangle ++(1,1);
    \path[fill=color-rf] (6,8) rectangle ++(1,1);
    \path[fill=color-rf] (7,5) rectangle ++(1,1);
    \path[fill=color-rf] (8,6) rectangle ++(1,1);
    \path[fill=color-r2f] (1,2) rectangle ++(1,1);
    \path[fill=color-r2f] (2,1) rectangle ++(1,1);
    \path[fill=color-r2f] (3,4) rectangle ++(1,1);
    \path[fill=color-r2f] (4,3) rectangle ++(1,1);
    \path[fill=color-r2f] (5,6) rectangle ++(1,1);
    \path[fill=color-r2f] (6,7) rectangle ++(1,1);
    \path[fill=color-r2f] (7,8) rectangle ++(1,1);
    \path[fill=color-r2f] (8,5) rectangle ++(1,1);
    \path[fill=color-r3f] (1,1) rectangle ++(1,1);
    \path[fill=color-r3f] (2,4) rectangle ++(1,1);
    \path[fill=color-r3f] (3,3) rectangle ++(1,1);
    \path[fill=color-r3f] (4,2) rectangle ++(1,1);
    \path[fill=color-r3f] (5,5) rectangle ++(1,1);
    \path[fill=color-r3f] (6,6) rectangle ++(1,1);
    \path[fill=color-r3f] (7,7) rectangle ++(1,1);
    \path[fill=color-r3f] (8,8) rectangle ++(1,1);
    \foreach \i in {1,...,\n} {
      \draw [very thin] (\i,1) -- (\i,\n); 
      \draw [very thin] (\i,\n+.1) -- (\i,\n+1.1); 
      \draw [very thin] (1,\i) -- (\n,\i); 
      \draw [very thin] (-.1,\i) -- (.9,\i); 
    }
    \draw [very thin] (1,\n+.1) rectangle (\n,\n+1.1);
    \draw [very thin] (-.1,1) rectangle (.9,\n);
    \node at (0.4,8.5) {$1$};
    \node at (0.4,7.5) {$r$};
    \node at (0.4,6.5) {$r^2$};
    \node at (0.4,5.5) {$r^3$}; 
    \node at (0.4,4.5) {$f$}; 
    \node at (0.4,3.5) {$rf$};
    \node at (0.4,2.5) {$r^2\!f$};
    \node at (0.4,1.5) {$r^3\!f$};
    \node at (1.5,9.6) {$1$};
    \node at (2.5,9.6) {$r$};
    \node at (3.5,9.6) {$r^2$};
    \node at (4.5,9.6) {$r^3$}; 
    \node at (5.5,9.6) {$f$}; 
    \node at (6.5,9.6) {$rf$};
    \node at (7.5,9.6) {$r^2\!f$};
    \node at (8.5,9.6) {$r^3\!f$};
    \node at (1.5,8.5) {$1$};
    \node at (1.5,7.5) {$r$};
    \node at (1.5,6.5) {$r^2$};
    \node at (1.5,5.5) {$r^3$}; 
    \node at (1.5,4.5) {$f$}; 
    \node at (1.5,3.5) {$rf$};
    \node at (1.5,2.5) {$r^2\!f$};
    \node at (1.5,1.5) {$r^3\!f$};
    \node at (2.5,8.5) {$r$};
    \node at (2.5,7.5) {$r^2$};
    \node at (2.5,6.5) {$r^3$};
    \node at (2.5,5.5) {$1$}; 
    \node at (2.5,4.5) {$r^3\!f$}; 
    \node at (2.5,3.5) {$f$};
    \node at (2.5,2.5) {$rf$};
    \node at (2.5,1.5) {$r^2\!f$};
    \node at (3.5,8.5) {$r^2$};
    \node at (3.5,7.5) {$r^3$};
    \node at (3.5,6.5) {$1$};
    \node at (3.5,5.5) {$r$}; 
    \node at (3.5,4.5) {$r^2\!f$}; 
    \node at (3.5,3.5) {$r^3\!f$};
    \node at (3.5,2.5) {$f$};
    \node at (3.5,1.5) {$rf$};
    \node at (4.5,8.5) {$r^3$};
    \node at (4.5,7.5) {$1$};
    \node at (4.5,6.5) {$r$};
    \node at (4.5,5.5) {$r^2$}; 
    \node at (4.5,4.5) {$rf$}; 
    \node at (4.5,3.5) {$r^2\!f$};
    \node at (4.5,2.5) {$r^3\!f$};
    \node at (4.5,1.5) {$f$};
    \node at (5.5,8.5) {$f$};
    \node at (5.5,7.5) {$rf$};
    \node at (5.5,6.5) {$r^2\!f$};
    \node at (5.5,5.5) {$r^3\!f$}; 
    \node at (5.5,4.5) {$1$}; 
    \node at (5.5,3.5) {$r$};
    \node at (5.5,2.5) {$r^2$};
    \node at (5.5,1.5) {$r^3$};
    \node at (6.5,8.5) {$rf$};
    \node at (6.5,7.5) {$r^2\!f$};
    \node at (6.5,6.5) {$r^3\!f$};
    \node at (6.5,5.5) {$f$}; 
    \node at (6.5,4.5) {$r^3$}; 
    \node at (6.5,3.5) {$1$};
    \node at (6.5,2.5) {$r$};
    \node at (6.5,1.5) {$r^2$};
    \node at (7.5,8.5) {$r^2\!f$};
    \node at (7.5,7.5) {$r^3\!f$};
    \node at (7.5,6.5) {$f$};
    \node at (7.5,5.5) {$rf$}; 
    \node at (7.5,4.5) {$r^2$}; 
    \node at (7.5,3.5) {$r^3$};
    \node at (7.5,2.5) {$1$};
    \node at (7.5,1.5) {$r$};
    \node at (8.5,8.5) {$r^3\!f$};
    \node at (8.5,7.5) {$f$};
    \node at (8.5,6.5) {$rf$};
    \node at (8.5,5.5) {$r^2\!f$}; 
    \node at (8.5,4.5) {$r$}; 
    \node at (8.5,3.5) {$r^2$};
    \node at (8.5,2.5) {$r^3$};
    \node at (8.5,1.5) {$1$};
    \node at (5,.3) {\normalsize $D_4$};
  \end{scope}
  \begin{scope}[shift={(10.5,0)},box/.style={anchor=south}]
      \newcommand*{\m}{5}
      \colorlet{color-e}{tGreen!80!black}
      \colorlet{color-v}{tRed}
      \colorlet{color-h}{tBlue}
      \colorlet{color-vh}{tYellow}
      \tikzstyle{every node}=[font=\footnotesize]
      \path[fill=color-e] (-.1,4) rectangle ++(1,1);
      \path[fill=color-v] (-.1,3) rectangle ++(1,1);
      \path[fill=color-h] (-.1,2) rectangle ++(1,1);
      \path[fill=color-vh] (-.1,1) rectangle ++(1,1);
      \path[fill=color-e] (1,5.1) rectangle ++(1,1);
      \path[fill=color-v] (2,5.1) rectangle ++(1,1);
      \path[fill=color-h] (3,5.1) rectangle ++(1,1);
      \path[fill=color-vh] (4,5.1) rectangle ++(1,1);
      \path[fill=color-e] (1,4) rectangle ++(1,1);
      \path[fill=color-e] (2,3) rectangle ++(1,1);
      \path[fill=color-e] (3,2) rectangle ++(1,1);
      \path[fill=color-e] (4,1) rectangle ++(1,1);
      \path[fill=color-v] (1,3) rectangle ++(1,1);
      \path[fill=color-v] (2,4) rectangle ++(1,1);
      \path[fill=color-v] (3,1) rectangle ++(1,1);
      \path[fill=color-v] (4,2) rectangle ++(1,1);
      \path[fill=color-h] (1,2) rectangle ++(1,1);
      \path[fill=color-h] (2,1) rectangle ++(1,1);
      \path[fill=color-h] (3,4) rectangle ++(1,1);
      \path[fill=color-h] (4,3) rectangle ++(1,1);
      \path[fill=color-vh] (1,1) rectangle ++(1,1);
      \path[fill=color-vh] (2,2) rectangle ++(1,1);
      \path[fill=color-vh] (3,3) rectangle ++(1,1);
      \path[fill=color-vh] (4,4) rectangle ++(1,1);
      \foreach \i in {1,...,\m} {
        \draw [very thin] (\i,1) -- (\i,\m); 
        \draw [very thin] (\i,\m+.1) -- (\i,\m+1.1); 
        \draw [very thin] (1,\i) -- (\m,\i); 
        \draw [very thin] (-.1,\i) -- (.9,\i); 
      } 
      \draw [very thin] (1,\m+.1) rectangle (\m,\m+1.1);
      \draw [very thin] (-.1,1) rectangle (.9,\m); 
      \node [box] at (0.4,4.1) {$\pm 1$}; 
      \node [box] at (0.4,3.1) {$\pm i$};
      \node [box] at (0.4,2.1) {$\pm j$};
      \node [box] at (0.4,1.1) {$\pm k$};
      \node [box] at (1.5,5.2) {$\pm 1$}; 
      \node [box] at (2.5,5.2) {$\pm i$};
      \node [box] at (3.5,5.2) {$\pm j$};
      \node [box] at (4.5,5.2) {$\pm k$};
      \node [box] at (1.5,4.1) {$\pm 1$}; 
      \node [box] at (1.5,3.1) {$\pm i$};
      \node [box] at (1.5,2.1) {$\pm j$};
      \node [box] at (1.5,1.1) {$\pm k$};
      \node [box] at (2.5,4.1) {$\pm i$}; 
      \node [box] at (2.5,3.1) {$\pm 1$};
      \node [box] at (2.5,2.1) {$\pm k$};
      \node [box] at (2.5,1.1) {$\pm j$};
      \node [box] at (3.5,4.1) {$\pm j$}; 
      \node [box] at (3.5,3.1) {$\pm k$};
      \node [box] at (3.5,2.1) {$\pm 1$};
      \node [box] at (3.5,1.1) {$\pm i$};
      \node [box] at (4.5,4.1) {$\pm k$}; 
      \node [box] at (4.5,3.1) {$\pm j$};
      \node [box] at (4.5,2.1) {$\pm i$};
      \node [box] at (4.5,1.1) {$\pm 1$};
      \node at (3,.3) {\normalsize $Q_8/\<-1\>\cong D_2$};
    \end{scope}
    \begin{scope}[shift={(17,0)}]
    \tikzstyle{every node}=[font=\scriptsize]
    \path[fill=color-1] (-.1,8) rectangle ++(1,1);
    \path[fill=color-r] (-.1,7) rectangle ++(1,1);
    \path[fill=color-r2] (-.1,6) rectangle ++(1,1);
    \path[fill=color-r3] (-.1,5) rectangle ++(1,1);
    \path[fill=color-f] (-.1,4) rectangle ++(1,1);
    \path[fill=color-rf] (-.1,3) rectangle ++(1,1);
    \path[fill=color-r2f] (-.1,2) rectangle ++(1,1);
    \path[fill=color-r3f] (-.1,1) rectangle ++(1,1);
    \path[fill=color-1] (1,9.1) rectangle ++(1,1);
    \path[fill=color-r] (2,9.1) rectangle ++(1,1);
    \path[fill=color-r2] (3,9.1) rectangle ++(1,1);
    \path[fill=color-r3] (4,9.1) rectangle ++(1,1);
    \path[fill=color-f] (5,9.1) rectangle ++(1,1);
    \path[fill=color-rf] (6,9.1) rectangle ++(1,1);
    \path[fill=color-r2f] (7,9.1) rectangle ++(1,1);
    \path[fill=color-r3f] (8,9.1) rectangle ++(1,1);
    \path[fill=color-1] (1,8) rectangle ++(1,1);
    \path[fill=color-1] (2,5) rectangle ++(1,1);
    \path[fill=color-1] (3,6) rectangle ++(1,1);
    \path[fill=color-1] (4,7) rectangle ++(1,1);
    \path[fill=color-1] (5,4) rectangle ++(1,1);
    \path[fill=color-1] (6,3) rectangle ++(1,1);
    \path[fill=color-1] (7,2) rectangle ++(1,1);
    \path[fill=color-1] (8,1) rectangle ++(1,1);
    \path[fill=color-r] (1,7) rectangle ++(1,1);
    \path[fill=color-r] (2,8) rectangle ++(1,1);
    \path[fill=color-r] (3,5) rectangle ++(1,1);
    \path[fill=color-r] (4,6) rectangle ++(1,1);
    \path[fill=color-r] (5,3) rectangle ++(1,1);
    \path[fill=color-r] (6,2) rectangle ++(1,1);
    \path[fill=color-r] (7,1) rectangle ++(1,1);
    \path[fill=color-r] (8,4) rectangle ++(1,1);
    \path[fill=color-r2] (1,6) rectangle ++(1,1);
    \path[fill=color-r2] (2,7) rectangle ++(1,1);
    \path[fill=color-r2] (3,8) rectangle ++(1,1);
    \path[fill=color-r2] (4,5) rectangle ++(1,1);
    \path[fill=color-r2] (5,2) rectangle ++(1,1);
    \path[fill=color-r2] (6,1) rectangle ++(1,1);
    \path[fill=color-r2] (7,4) rectangle ++(1,1);
    \path[fill=color-r2] (8,3) rectangle ++(1,1);
    \path[fill=color-r3] (1,5) rectangle ++(1,1);
    \path[fill=color-r3] (2,6) rectangle ++(1,1);
    \path[fill=color-r3] (3,7) rectangle ++(1,1);
    \path[fill=color-r3] (4,8) rectangle ++(1,1);
    \path[fill=color-r3] (5,1) rectangle ++(1,1);
    \path[fill=color-r3] (6,4) rectangle ++(1,1);
    \path[fill=color-r3] (7,3) rectangle ++(1,1);
    \path[fill=color-r3] (8,2) rectangle ++(1,1);
    \path[fill=color-f] (1,4) rectangle ++(1,1);
    \path[fill=color-f] (2,3) rectangle ++(1,1);
    \path[fill=color-f] (3,2) rectangle ++(1,1);
    \path[fill=color-f] (4,1) rectangle ++(1,1);
    \path[fill=color-f] (5,8) rectangle ++(1,1);
    \path[fill=color-f] (6,5) rectangle ++(1,1);
    \path[fill=color-f] (7,6) rectangle ++(1,1);
    \path[fill=color-f] (8,7) rectangle ++(1,1);
    \path[fill=color-rf] (1,3) rectangle ++(1,1);
    \path[fill=color-rf] (2,2) rectangle ++(1,1);
    \path[fill=color-rf] (3,1) rectangle ++(1,1);
    \path[fill=color-rf] (4,4) rectangle ++(1,1);
    \path[fill=color-rf] (5,7) rectangle ++(1,1);
    \path[fill=color-rf] (6,8) rectangle ++(1,1);
    \path[fill=color-rf] (7,5) rectangle ++(1,1);
    \path[fill=color-rf] (8,6) rectangle ++(1,1);
    \path[fill=color-r2f] (1,2) rectangle ++(1,1);
    \path[fill=color-r2f] (2,1) rectangle ++(1,1);
    \path[fill=color-r2f] (3,4) rectangle ++(1,1);
    \path[fill=color-r2f] (4,3) rectangle ++(1,1);
    \path[fill=color-r2f] (5,6) rectangle ++(1,1);
    \path[fill=color-r2f] (6,7) rectangle ++(1,1);
    \path[fill=color-r2f] (7,8) rectangle ++(1,1);
    \path[fill=color-r2f] (8,5) rectangle ++(1,1);
    \path[fill=color-r3f] (1,1) rectangle ++(1,1);
    \path[fill=color-r3f] (2,4) rectangle ++(1,1);
    \path[fill=color-r3f] (3,3) rectangle ++(1,1);
    \path[fill=color-r3f] (4,2) rectangle ++(1,1);
    \path[fill=color-r3f] (5,5) rectangle ++(1,1);
    \path[fill=color-r3f] (6,6) rectangle ++(1,1);
    \path[fill=color-r3f] (7,7) rectangle ++(1,1);
    \path[fill=color-r3f] (8,8) rectangle ++(1,1);
    \foreach \i in {1,...,\n} {
      \draw [very thin] (\i,1) -- (\i,\n); 
      \draw [very thin] (\i,\n+.1) -- (\i,\n+1.1); 
      \draw [very thin] (1,\i) -- (\n,\i); 
      \draw [very thin] (-.1,\i) -- (.9,\i); 
    }
    \draw [very thin] (1,\n+.1) rectangle (\n,\n+1.1);
    \draw [very thin] (-.1,1) rectangle (.9,\n);
    \node at (0.4,8.5) {$\pm 1$};
    \node at (0.4,7.5) {$\pm\zeta$};
    \node at (0.4,6.5) {$\pm\zeta^2$};
    \node at (0.4,5.5) {$\pm\zeta^3$}; 
    \node at (0.4,4.5) {$j$}; 
    \node at (0.4,3.5) {$\pm\zeta j$};
    \node at (0.4,2.5) {$\pm\zeta^2\!j$};
    \node at (0.4,1.5) {$\pm\zeta^3\!j$};
    \node at (1.5,9.6) {$\pm 1$};
    \node at (2.5,9.6) {$\pm\zeta$};
    \node at (3.5,9.6) {$\pm\zeta^2$};
    \node at (4.5,9.6) {$\pm\zeta^3$}; 
    \node at (5.5,9.6) {$j$}; 
    \node at (6.5,9.6) {$\pm\zeta j$};
    \node at (7.5,9.6) {$\pm\zeta^2\!j$};
    \node at (8.5,9.6) {$\pm\zeta^3\!j$};
    \node at (1.5,8.5) {$\pm 1$};
    \node at (1.5,7.5) {$\pm\zeta$};
    \node at (1.5,6.5) {$\pm\zeta^2$};
    \node at (1.5,5.5) {$\pm\zeta^3$}; 
    \node at (1.5,4.5) {$j$}; 
    \node at (1.5,3.5) {$\pm\zeta j$};
    \node at (1.5,2.5) {$\pm\zeta^2\!j$};
    \node at (1.5,1.5) {$\pm\zeta^3\!j$};
    \node at (2.5,8.5) {$\pm\zeta$};
    \node at (2.5,7.5) {$\pm\zeta^2$};
    \node at (2.5,6.5) {$\pm\zeta^3$};
    \node at (2.5,5.5) {$\pm 1$}; 
    \node at (2.5,4.5) {$\pm\zeta^3\!j$}; 
    \node at (2.5,3.5) {$j$};
    \node at (2.5,2.5) {$\pm\zeta j$};
    \node at (2.5,1.5) {$\pm\zeta^2\!j$};
    \node at (3.5,8.5) {$\pm\zeta^2$};
    \node at (3.5,7.5) {$\pm\zeta^3$};
    \node at (3.5,6.5) {$\pm 1$};
    \node at (3.5,5.5) {$\pm\zeta$}; 
    \node at (3.5,4.5) {$\pm\zeta^2\!j$}; 
    \node at (3.5,3.5) {$\pm\zeta^3\!j$};
    \node at (3.5,2.5) {$j$};
    \node at (3.5,1.5) {$\pm\zeta j$};
    \node at (4.5,8.5) {$\pm\zeta^3$};
    \node at (4.5,7.5) {$\pm 1$};
    \node at (4.5,6.5) {$\pm\zeta$};
    \node at (4.5,5.5) {$\pm\zeta^2$}; 
    \node at (4.5,4.5) {$\pm\zeta j$}; 
    \node at (4.5,3.5) {$\pm\zeta^2\!j$};
    \node at (4.5,2.5) {$\pm\zeta^3\!j$};
    \node at (4.5,1.5) {$j$};
    \node at (5.5,8.5) {$j$};
    \node at (5.5,7.5) {$\pm\zeta j$};
    \node at (5.5,6.5) {$\pm\zeta^2\!j$};
    \node at (5.5,5.5) {$\pm\zeta^3\!j$}; 
    \node at (5.5,4.5) {$\pm 1$}; 
    \node at (5.5,3.5) {$\pm\zeta$};
    \node at (5.5,2.5) {$\pm\zeta^2$};
    \node at (5.5,1.5) {$\pm\zeta^3$};
    \node at (6.5,8.5) {$\pm\zeta j$};
    \node at (6.5,7.5) {$\pm\zeta^2\!j$};
    \node at (6.5,6.5) {$\pm\zeta^3\!j$};
    \node at (6.5,5.5) {$j$}; 
    \node at (6.5,4.5) {$\pm\zeta^3$}; 
    \node at (6.5,3.5) {$\pm 1$};
    \node at (6.5,2.5) {$\pm\zeta$};
    \node at (6.5,1.5) {$\pm\zeta^2$};
    \node at (7.5,8.5) {$\pm\zeta^2\!j$};
    \node at (7.5,7.5) {$\pm\zeta^3\!j$};
    \node at (7.5,6.5) {$j$};
    \node at (7.5,5.5) {$\pm\zeta j$}; 
    \node at (7.5,4.5) {$\pm\zeta^2$}; 
    \node at (7.5,3.5) {$\pm\zeta^3$};
    \node at (7.5,2.5) {$\pm 1$};
    \node at (7.5,1.5) {$\pm\zeta$};
    \node at (8.5,8.5) {$\pm\zeta^3\!j$};
    \node at (8.5,7.5) {$j$};
    \node at (8.5,6.5) {$\pm\zeta j$};
    \node at (8.5,5.5) {$\pm\zeta^2\!j$}; 
    \node at (8.5,4.5) {$\pm\zeta$}; 
    \node at (8.5,3.5) {$\pm\zeta^2$};
    \node at (8.5,2.5) {$\pm\zeta^3$};
    \node at (8.5,1.5) {$\pm 1$};
    \node at (5,.3) {\normalsize $Q_{16}/\<-1\>\cong D_4$};
  \end{scope}  
  \end{tikzpicture}
  \caption{Left: A Cayley table of the dihedral group. Middle: the quotient $Q_8/\<-1\>$ is isomorphic to the Klein $4$-group, which is also the dihedral group $D_2$. Right: More generally, $Q_{2^n}/\<-1\>$ is dihedral.}\label{fig:cayley-tables}
  \end{figure}

Another useful visual aid for groups that is typically underutilized in algebra classes are \emph{Cayley tables}. These are basically multiplication tables but for groups, and should be self-explanatory. However, they can reveal patterns that Cayley graphs might miss. On the left in Figure~\ref{fig:cayley-tables} is the Cayley table for the dihedral group $D_4$. In the middle is a Cayley table for the quotient of $Q_8$ by the normal subgroup $\<-1\>$. The entries in this table are cosets of $\<-1\>=\{1,-1\}$, and it highlights how the quotient $Q_8/\<-1\>$ is the Klein $4$-group, also known as the dihedral group $D_2$. On the right in Figure~\ref{fig:cayley-tables} is the analogous quotient of $Q_{16}$ by $\<-1\>$. Notice how this has the structure of $D_4$, and so $Q_{16}/\<-1\>\cong D_4$. More generally, $Q_{2^n}/\<-1\>\cong D_{2^{n-2}}$.

%% Cayley tables for D_4, showing a quotient
%%

Though Cayley tables have their limitations, especially because they are simply not practical for large groups, the patterns within them have inspired various works of art. For example, Cayley graph quilts have been made, where quotient groups are visually apparent; see~\cite{jensen2023group}. Even though Cayley tables explicitly specify how to multiply two elements, Cayley graphs also serve as a ``group calculator,'' making it easy to compute the product of any two elements. They also highlight patterns and symmetries that are not necessarily apparent from just a presentation, despite these being no more than an algebraic encoding of a Cayley graph, albeit without the information of the size of the group.

In addition to the aesthetics of Cayley graphs, another advantage of them is that they elucidate groups whose structure is muddied by a complicated presentation. One fantastic example of this is the following group of order $16$, defined by the presentation
\[
G=\big\<a,b,c\mid a^4=c^2=1,\,a^2=b^2,\; ab=ba,\,ac=ca,\,a^2b=cbc\big\>.
\]
This group is sometimes called the \emph{Pauli group} (on $1$ qubit), because it can be generated by the so-called \emph{Pauli matrices} from quantum physics. Not surprisingly, it gets little-to-no attention in a standard algebra class or book. Nick Webb, a retired scientist and recreational mathematician who maintains a fantastic website called \emph{Weddslist}~\cite{wedd2024weddslist}, containing many beautiful pictures of groups and Cayley graphs, writes:
\begin{quote}
``\emph{Among the groups of order less than 32, the Pauli group is one of the hardest to understand\dots it has $C_4\times C_2$ and $D_4$ and $Q_8$ all as normal subgroups. I don't know about my readers, but I find that hard to envisage. This page aims to make the structure of the Pauli group easier to understand.}'' 
\end{quote}

\begin{figure}[!ht]
  \begin{tikzpicture}[scale=1.6,auto]
  \tikzstyle{v} = [circle, draw, fill=lightgray,inner sep=0pt, 
  minimum size=7mm]
\tikzstyle{v-yel} = [circle, draw, fill=vYellow,inner sep=0pt,
  minimum size=7mm]
\tikzstyle{rr-bend} = [draw, very thick, eRed,bend left=12]
\tikzstyle{bb-bend} = [draw, very thick, eBlue,bend left=12]
    \tikzstyle{every node}=[font=\tiny]
    \begin{scope}[shift={(0,0)}]
    \setlength{\arraycolsep}{1.5pt}
    \renewcommand{\arraystretch}{.7}
      \node (s) at (90:1.1) [v] {$\begin{bmatrix}i&0\\0&i\end{bmatrix}$};
      \node (rs) at (45:1.15) [v-yel] {$\begin{bmatrix}0&-\!1\\1&0\end{bmatrix}$};
      \node (r2s) at (0:1.15) [v] {$\begin{bmatrix}1&0\\0&-\!1\end{bmatrix}$};
      \node (r3s) at (-45:1.15) [v-yel] {$\begin{bmatrix}0&i\\i&0\end{bmatrix}$};
      \node (r4s) at (-90:1.15) [v] {$\begin{bmatrix}-\!i&0\\0&\!-\!i\end{bmatrix}$};
      \node (r5s) at (-135:1.15) [v-yel] {$\begin{bmatrix}0&1\\-\!1&0\end{bmatrix}$};
      \node (r6s) at (180:1.15) [v] {$\begin{bmatrix}-\!1&0\\0&1\end{bmatrix}$};
      \node (r7s) at (135:1.15) [v-yel] {$\begin{bmatrix}0&-\!i\\\!\!-\!i&0\end{bmatrix}$};
      \node (1) at (90:2) [v-yel] {$\begin{bmatrix}i&0\\0&-\!i\end{bmatrix}$};
      \node (r) at (45:2) [v] {$\begin{bmatrix}0&1\\1&0\end{bmatrix}$};
      \node (r2) at (0:2) [v-yel] {$\begin{bmatrix}1&0\\0&1\end{bmatrix}$};
      \node (r3) at (-45:2) [v] {$\begin{bmatrix}0&-\!i\\i&0\end{bmatrix}$};
      \node (r4) at (-90:2) [v-yel] {$\begin{bmatrix}-\!i&0\\0&i\end{bmatrix}$};
      \node (r5) at (-135:2) [v] {\setlength{\arraycolsep}{1.1pt}$\begin{bmatrix}0&\!\!-\!1\\-\!1&0\end{bmatrix}$};
      \node (r6) at (180:2) [v-yel] {\setlength{\arraycolsep}{1.1pt}$\begin{bmatrix}-\!1&0\\0&\!\!-\!1\end{bmatrix}$};
      \node (r7) at (135:2) [v] {$\begin{bmatrix}0&i\\-\!i&0\end{bmatrix}$};
      \node at (9:2.28) {\large\Palert{$\mathbf{1}$}};
      \node at (90+8:2.3) {\large\Palert{$\mathbf{i}$}};
      \node at (90+30:1.4) {\large\Palert{$\mathbf{k}$}};
      \node at (90-30:1.4) {\large\Palert{$\mathbf{j}$}};
      \node at (180-9:2.35) {\large\Palert{$\mathbf{-1}$}};
      \node at (270-8:2.35) {\large\Palert{$\mathbf{-i}$}};
      \node at (270-30:1.4) {\large\Palert{$\mathbf{-j}$}};
      \node at (270+30:1.4) {\large\Palert{$\mathbf{-k}$}};
      \draw [rr-bend] (1) to (r);
      \draw [bb-bend] (r) to (r2);
      \draw [rr-bend] (r2) to (r3);
      \draw [bb-bend] (r3) to (r4);
      \draw [rr-bend] (r4) to (r5);
      \draw [bb-bend] (r5) to (r6);
      \draw [rr-bend] (r6) to (r7);
      \draw [bb-bend] (r7) to (1);
      \draw [bb] (s) to (r3s);
      \draw [rr] (r3s) to (r6s);
      \draw [bb] (r6s) to (rs);
      \draw [rr] (rs) to (r4s);
      \draw [bb] (r4s) to (r7s);
      \draw [rr] (r7s) to (r2s);
      \draw [bb] (r2s) to (r5s);
      \draw [rr] (r5s) to (s);
      \draw [gg] (1) to (s); \draw [gg] (r) to (rs);
      \draw [gg] (r2) to (r2s); \draw [gg] (r3) to (r3s);
      \draw [gg] (r4) to (r4s); \draw [gg] (r5) to (r5s);
      \draw [gg] (r6) to (r6s); \draw [gg] (r7) to (r7s);
      \node at (0,0) {\normalsize $\DQ_8$};
    \end{scope}
      \begin{scope}[shift={(5,0)},scale=1.05]
      \tikzstyle{v} = [circle, draw, fill=lightgray,inner sep=0pt, 
  minimum size=4mm]
\tikzstyle{v-yel} = [circle, draw, fill=vYellow,inner sep=0pt, 
  minimum size=4mm]
\tikzstyle{rr-bend} = [draw, very thick, eRed,bend left=7]
\tikzstyle{bb-bend} = [draw, very thick, eBlue,bend left=7]
\tikzstyle{every node}=[font=\scriptsize]
      %\node at (45:2.3) {\large\Palert{$\zeta_8$}};
      %\node at (30:1.45) {\large\Palert{$j$}};
      \node (s) at (90:1.2) [v] {};
      \node (rs) at (67.5:1.2) [v-yel] {};
      \node (r2s) at (45:1.2) [v] {};
      \node (r3s) at (22.5:1.2) [v-yel] {$j$};
      \node (r4s) at (0:1.2) [v] {};
      \node (r5s) at (-22.5:1.2) [v-yel] {};
      \node (r6s) at (-45:1.2) [v] {};
      \node (r7s) at (-67.5:1.2) [v-yel] {};
      \node (r8s) at (-90:1.2) [v] {};
      \node (r9s) at (-112.5:1.2) [v-yel] {};
      \node (r10s) at (-135:1.2) [v] {};
      \node (r11s) at (-157.5:1.2) [v-yel] {};
      \node (r12s) at (180:1.2) [v] {};
      \node (r13s) at (157.5:1.2) [v-yel] {};
      \node (r14s) at (135:1.2) [v] {};
      \node (r15s) at (112.5:1.2) [v-yel] {};
      \node (1) at (90:2) [v-yel] {$i$};
      \node (r) at (67.5:2) [v] {};
      \node (r2) at (45:2) [v-yel] {$\zeta$};
      \node (r3) at (22.5:2) [v] {$f$};
      \node (r4) at (0:2) [v-yel] {$I$};
      \node (r5) at (-22.5:2) [v] {};
      \node (r6) at (-45:2) [v-yel] {};
      \node (r7) at (-67.5:2) [v] {};
      \node (r8) at (-90:2) [v-yel] {};
      \node (r9) at (-112.5:2) [v] {};
      \node (r10) at (-135:2) [v-yel] {};
      \node (r11) at (-157.5:2) [v] {};
      \node (r12) at (180:2) [v-yel] {};
      \node (r13) at (157.5:2) [v] {};
      \node (r14) at (135:2) [v-yel] {};
      \node (r15) at (112.5:2) [v] {};
      \draw [rr-bend] (1) to (r); \draw [bb-bend] (r) to (r2);
      \draw [rr-bend] (r2) to (r3); \draw [bb-bend] (r3) to (r4);
      \draw [rr-bend] (r4) to (r5); \draw [bb-bend] (r5) to (r6);
      \draw [rr-bend] (r6) to (r7); \draw [bb-bend] (r7) to (r8);
      \draw [rr-bend] (r8) to (r9); \draw [bb-bend] (r9) to (r10);
      \draw [rr-bend] (r10) to (r11); \draw [bb-bend] (r11) to (r12);
      \draw [rr-bend] (r12) to (r13); \draw [bb-bend] (r13) to (r14);
      \draw [rr-bend] (r14) to (r15); \draw [bb-bend] (r15) to (1);
      \draw [rr] (s) to (r9s); \draw [bb] (r9s) to (r2s);
      \draw [rr] (r2s) to (r11s); \draw [bb] (r11s) to (r4s);
      \draw [rr] (r4s) to (r13s); \draw [bb] (r13s) to (r6s);
      \draw [rr] (r6s) to (r15s); \draw [bb] (r15s) to (r8s);
      \draw [rr] (r8s) to (rs); \draw [bb] (rs) to (r10s);
      \draw [rr] (r10s) to (r3s); \draw [bb] (r3s) to (r12s);
      \draw [rr] (r12s) to (r5s); \draw [bb] (r5s) to (r14s);
      \draw [rr] (r14s) to (r7s); \draw [bb] (r7s) to (s); 
      \draw [gg] (1) to (s); \draw [gg] (r) to (rs);
      \draw [gg] (r2) to (r2s); \draw [gg] (r3) to (r3s);
      \draw [gg] (r4) to (r4s); \draw [gg] (r5) to (r5s);
      \draw [gg] (r6) to (r6s); \draw [gg] (r7) to (r7s);
      \draw [gg] (r8) to (r8s); \draw [gg] (r9) to (r9s);
      \draw [gg] (r10) to (r10s); \draw [gg] (r11) to (r11s);
      \draw [gg] (r12) to (r12s); \draw [gg] (r13) to (r13s);
      \draw [gg] (r14) to (r14s); \draw [gg] (r15) to (r15s);
      \node at (0,0) {\scriptsize $\DQ_{16}$};
    \end{scope}
  \end{tikzpicture}
\caption{To construct the diquaternion group $\DQ_8=\<i,j,f\>$, take the standard representation of $Q_8$, and throw in the standard ``reflection matrix'' from $D_n$.}\label{fig:DQ8-DQ16}
\end{figure}
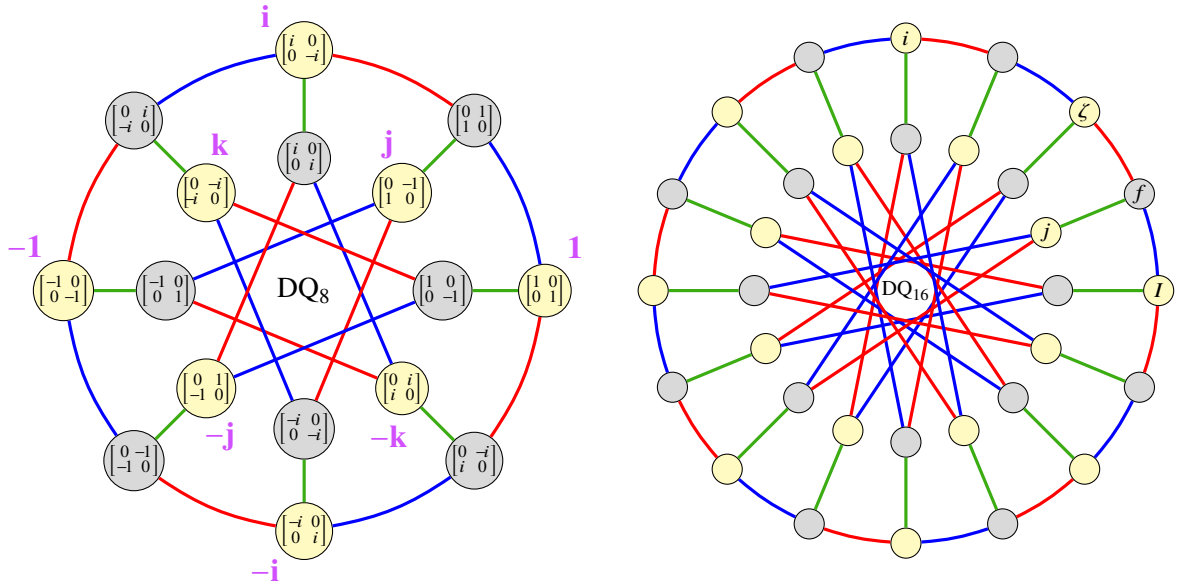

On his webpage, Wedd includes a number of aesthetically pleasing pictures of different Cayley graphs of the Pauli group, including several embedded on a torus, and another on a double-torus. An alternate construction of this group involves creatively ``combining'' the dihedral and quaternion groups from Figure~\ref{fig:D4-Q8-Q16}. Specifically, these groups can be represented by $2\times 2$ complex matrices, via the following maps:
\[
r\mapsto\begin{bmatrix}e^{2\pi i/4}&0\\0&e^{-2\pi i/4}\end{bmatrix},\qquad
i\mapsto\begin{bmatrix}i&0\\0&-i\end{bmatrix},\qquad
j\mapsto\begin{bmatrix}0&-1\\1&0\end{bmatrix},\qquad
f\mapsto\begin{bmatrix}0&1\\1&0\end{bmatrix}.
\]
Since $i=e^{2\pi i/4}$, the matrices for $r$ and $i$ coincide. The Pauli group is the group generated by the three matrices above, a process that can be playfully called \emph{dihedralizing the quaternions}~\cite{macauley2024dihedralizing}. This also gives a better name to this group: the \emph{diquaternion group}, denoted $\DQ_8$. This name has also never appeared in the literature, except in the citation above, which is a paper that is in press at the time of the submission of this paper. A Cayley graph of $\DQ_8$ is shown on the left in Figure~\ref{fig:DQ8-DQ16}. On the right is the result of dihedralizing the (generalized) quaternion group $Q_{16}$, by replacing the $e^{2\pi i/4}$ entries in the matrix above with $\zeta_8=e^{2\pi i/8}$. A Cayley graph is shown in Figure~\ref{fig:DQ8-DQ16} on the right.

%% DQ_8 (Pauli group) Cayley graph
%%

Upon seeing the Cayley graph, the three aforementioned index-$2$ subgroups in $\DQ_8$ can be seen by inspection: \emph{(i)} the matrices in $Q_8$ are highlighted \emph{(ii)} the matrices in the ``outer ring'' comprise $D_4$, and \emph{(iii)} the matrices ``along the $x$- and $y$-axes'' define $C_4\times C_2$. Similarly, the group $\DQ_{16}$ has analogous index-$2$ subgroups $Q_{16}$, $D_8$, and $C_8\times C_2$. 

Returning to the quaternion group $Q_8$ and generalized quaternion group $Q_{16}$ in Figure~\ref{fig:D4-Q8-Q16}, it should be apparent how to define a generalized quaternion group $Q_{2^n}$ for any $n\geq 4$. Note that $\<r\>$ is always an index-$2$ cyclic subgroup, but $Q_{2^n}$ is never a semidirect product of this with $C_2$. In fact,  it is straightforward to show that every cyclic subgroup $\<g\>$ contains $-1$, and so $Q_{2^n}$ can never decompose as a semidirect product of \emph{any} of two of its nontrivial subgroups. 

It is worth mentioning that for each $n$, the \emph{Pauli group on $n$ qubits} can be constructed by taking tensors of the Pauli matrices. However, this has order $4^{n+1}$, and so the Pauli groups on $2$, $3$, and $4$ qubits have order $64$, $256$, and $1024$, respectively. In contrast, there is a diquaternion group of order $2^n$ for all $n\geq 4$.

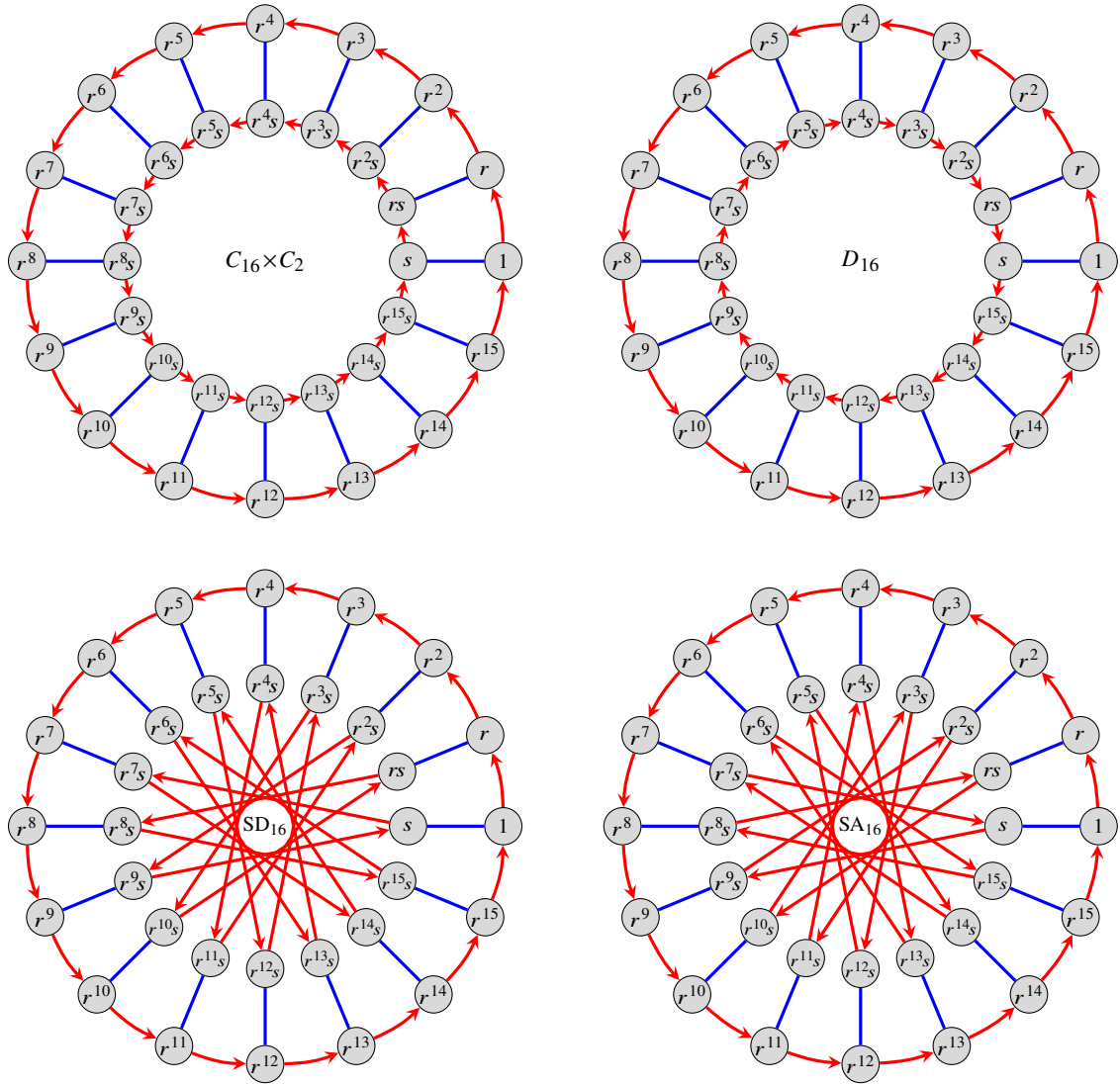
\begin{figure}[!ht]
  \begin{tikzpicture}[scale=1.6,auto]
\tikzstyle{v} = [circle, draw, fill=lightgray,inner sep=0pt, 
  minimum size=5mm]
  \tikzstyle{R} = [draw, very thick, eRed,-stealth,bend right=10]
\tikzstyle{R2} = [draw, very thick, eRed,-stealth,bend right=15]
    \tikzstyle{every node}=[font=\scriptsize]
    \begin{scope}[shift={(0,4.75)}]
      \node (1) at (0:2) [v] {$1$};
      \node (r) at (22.5:2) [v] {$r$};
      \node (r2) at (45:2) [v] {$r^2$};
      \node (r3) at (67.5:2) [v] {$r^3$};
      \node (r4) at (90:2) [v] {$r^4$};
      \node (r5) at (112.5:2) [v] {$r^5$};
      \node (r6) at (135:2) [v] {$r^6$};
      \node (r7) at (157.5:2) [v] {$r^7$};
      \node (r8) at (180:2) [v] {$r^8$};
      \node (r9) at (-157.5:2) [v] {$r^9$};
      \node (r10) at (-135:2) [v] {$r^{10}$};
      \node (r11) at (-112.5:2) [v] {$r^{11}$};
      \node (r12) at (-90:2) [v] {$r^{12}$};
      \node (r13) at (-67.5:2) [v] {$r^{13}$};
      \node (r14) at (-45:2) [v] {$r^{14}$};
      \node (r15) at (-22.5:2) [v] {$r^{15}$};
      \node (s) at (0:1.2) [v] {$s$};
      \node (rs) at (22.5:1.2) [v] {$r\!s$};
      \node (r2s) at (45:1.2) [v] {$r^2\!s$};
      \node (r3s) at (67.5:1.2) [v] {$r^3\!s$};
      \node (r4s) at (90:1.2) [v] {$r^4\!s$};
      \node (r5s) at (112.5:1.2) [v] {$r^5\!s$};
      \node (r6s) at (135:1.2) [v] {$r^6\!s$};
      \node (r7s) at (157.5:1.2) [v] {$r^7\!s$};
      \node (r8s) at (180:1.2) [v] {$r^8\!s$};
      \node (r9s) at (-157.5:1.2) [v] {$r^9\!s$};
       \tikzstyle{every node}=[font=\tiny]
      \node (r10s) at (-135:1.2) [v] {$r^{10}\!s$};
      \node (r11s) at (-112.5:1.2) [v] {$r^{11}\!s$};
      \node (r12s) at (-90:1.2) [v] {$r^{12}\!s$};
      \node (r13s) at (-67.5:1.2) [v] {$r^{13}\!s$};
      \node (r14s) at (-45:1.2) [v] {$r^{14}\!s$};
      \node (r15s) at (-22.5:1.2) [v] {$r^{15}\!s$};
      \draw [R] (1) to (r); \draw [R] (r) to (r2); \draw [R] (r2) to (r3);
      \draw [R] (r3) to (r4); \draw [R] (r4) to (r5); \draw [R] (r5) to (r6);
      \draw [R] (r6) to (r7); \draw [R] (r7) to (r8); \draw [R] (r8) to (r9);
      \draw [R] (r9) to (r10); \draw [R] (r10) to (r11);
      \draw [R] (r11) to (r12); \draw [R] (r12) to (r13);
      \draw [R] (r13) to (r14); \draw [R] (r14) to (r15);
      \draw [R] (r15) to (1);
      \draw [r] (s) to (rs); \draw [r] (rs) to (r2s);
      \draw [r] (r2s) to (r3s); \draw [r] (r3s) to (r4s);
      \draw [r] (r4s) to (r5s); \draw [r] (r5s) to (r6s);
      \draw [r] (r6s) to (r7s); \draw [r] (r7s) to (r8s);
      \draw [r] (r8s) to (r9s); \draw [r] (r9s) to (r10s);
      \draw [r] (r10s) to (r11s); \draw [r] (r11s) to (r12s);
      \draw [r] (r12s) to (r13s); \draw [r] (r13s) to (r14s);
      \draw [r] (r14s) to (r15s); \draw [r] (r15s) to (s); 
      \draw [bb] (1) to (s); \draw [bb] (r) to (rs);
      \draw [bb] (r2) to (r2s); \draw [bb] (r3) to (r3s);
      \draw [bb] (r4) to (r4s); \draw [bb] (r5) to (r5s);
      \draw [bb] (r6) to (r6s); \draw [bb] (r7) to (r7s);
      \draw [bb] (r8) to (r8s); \draw [bb] (r9) to (r9s);
      \draw [bb] (r10) to (r10s); \draw [bb] (r11) to (r11s);
      \draw [bb] (r12) to (r12s); \draw [bb] (r13) to (r13s);
      \draw [bb] (r14) to (r14s); \draw [bb] (r15) to (r15s);
      \node at (0,0) {\footnotesize $C_{16}\!\times\!C_2$};
    \end{scope}
    \begin{scope}[shift={(5,4.75)}]
      \tikzstyle{r-bend} = [draw, very thick, eRed,-stealth,bend right=8]
      \node (1) at (0:2) [v] {$1$};
      \node (r) at (22.5:2) [v] {$r$};
      \node (r2) at (45:2) [v] {$r^2$};
      \node (r3) at (67.5:2) [v] {$r^3$};
      \node (r4) at (90:2) [v] {$r^4$};
      \node (r5) at (112.5:2) [v] {$r^5$};
      \node (r6) at (135:2) [v] {$r^6$};
      \node (r7) at (157.5:2) [v] {$r^7$};
      \node (r8) at (180:2) [v] {$r^8$};
      \node (r9) at (-157.5:2) [v] {$r^9$};
      \node (r10) at (-135:2) [v] {$r^{10}$};
      \node (r11) at (-112.5:2) [v] {$r^{11}$};
      \node (r12) at (-90:2) [v] {$r^{12}$};
      \node (r13) at (-67.5:2) [v] {$r^{13}$};
      \node (r14) at (-45:2) [v] {$r^{14}$};
      \node (r15) at (-22.5:2) [v] {$r^{15}$};
      \node (s) at (0:1.2) [v] {$s$};
      \node (rs) at (22.5:1.2) [v] {$r\!s$};
      \node (r2s) at (45:1.2) [v] {$r^2\!s$};
      \node (r3s) at (67.5:1.2) [v] {$r^3\!s$};
      \node (r4s) at (90:1.2) [v] {$r^4\!s$};
      \node (r5s) at (112.5:1.2) [v] {$r^5\!s$};
      \node (r6s) at (135:1.2) [v] {$r^6\!s$};
      \node (r7s) at (157.5:1.2) [v] {$r^7\!s$};
      \node (r8s) at (180:1.2) [v] {$r^8\!s$};
      \node (r9s) at (-157.5:1.2) [v] {$r^9\!s$};
      \tikzstyle{every node}=[font=\tiny]
      \node (r10s) at (-135:1.2) [v] {$r^{10}\!s$};
      \node (r11s) at (-112.5:1.2) [v] {$r^{11}\!s$};
      \node (r12s) at (-90:1.2) [v] {$r^{12}\!s$};
      \node (r13s) at (-67.5:1.2) [v] {$r^{13}\!s$};
      \node (r14s) at (-45:1.2) [v] {$r^{14}\!s$};
      \node (r15s) at (-22.5:1.2) [v] {$r^{15}\!s$};

      \draw [r-bend] (1) to (r); \draw [r-bend] (r) to (r2);
      \draw [r-bend] (r2) to (r3); \draw [r-bend] (r3) to (r4);
      \draw [r-bend] (r4) to (r5); \draw [r-bend] (r5) to (r6);
      \draw [r-bend] (r6) to (r7); \draw [r-bend] (r7) to (r8);
      \draw [r-bend] (r8) to (r9); \draw [r-bend] (r9) to (r10);
      \draw [r-bend] (r10) to (r11); \draw [r-bend] (r11) to (r12);
      \draw [r-bend] (r12) to (r13); \draw [r-bend] (r13) to (r14);
      \draw [r-bend] (r14) to (r15); \draw [r-bend] (r15) to (1);
      \draw [r] (s) to (r15s); \draw [r] (r15s) to (r14s);
      \draw [r] (r14s) to (r13s); \draw [r] (r13s) to (r12s);
      \draw [r] (r12s) to (r11s); \draw [r] (r11s) to (r10s);
      \draw [r] (r10s) to (r9s); \draw [r] (r9s) to (r8s);
      \draw [r] (r8s) to (r7s); \draw [r] (r7s) to (r6s);
      \draw [r] (r6s) to (r5s); \draw [r] (r5s) to (r4s);
      \draw [r] (r4s) to (r3s); \draw [r] (r3s) to (r2s);
      \draw [r] (r2s) to (rs); \draw [r] (rs) to (s); 
      \draw [bb] (1) to (s); \draw [bb] (r) to (rs);
      \draw [bb] (r2) to (r2s); \draw [bb] (r3) to (r3s);
      \draw [bb] (r4) to (r4s); \draw [bb] (r5) to (r5s);
      \draw [bb] (r6) to (r6s); \draw [bb] (r7) to (r7s);
      \draw [bb] (r8) to (r8s); \draw [bb] (r9) to (r9s);
      \draw [bb] (r10) to (r10s); \draw [bb] (r11) to (r11s);
      \draw [bb] (r12) to (r12s); \draw [bb] (r13) to (r13s);
      \draw [bb] (r14) to (r14s); \draw [bb] (r15) to (r15s);
      \node at (0,0) {\footnotesize $D_{16}$};
    \end{scope}
    \begin{scope}[shift={(0,0)}]
      \node (1) at (0:2) [v] {$1$};
      \node (r) at (22.5:2) [v] {$r$};
      \node (r2) at (45:2) [v] {$r^2$};
      \node (r3) at (67.5:2) [v] {$r^3$};
      \node (r4) at (90:2) [v] {$r^4$};
      \node (r5) at (112.5:2) [v] {$r^5$};
      \node (r6) at (135:2) [v] {$r^6$};
      \node (r7) at (157.5:2) [v] {$r^7$};
      \node (r8) at (180:2) [v] {$r^8$};
      \node (r9) at (-157.5:2) [v] {$r^9$};
      \node (r10) at (-135:2) [v] {$r^{10}$};
      \node (r11) at (-112.5:2) [v] {$r^{11}$};
      \node (r12) at (-90:2) [v] {$r^{12}$};
      \node (r13) at (-67.5:2) [v] {$r^{13}$};
      \node (r14) at (-45:2) [v] {$r^{14}$};
      \node (r15) at (-22.5:2) [v] {$r^{15}$};
      \node (s) at (0:1.2) [v] {$s$};
      \node (rs) at (22.5:1.2) [v] {$r\!s$};
      \node (r2s) at (45:1.2) [v] {$r^2\!s$};
      \node (r3s) at (67.5:1.2) [v] {$r^3\!s$};
      \node (r4s) at (90:1.2) [v] {$r^4\!s$};
      \node (r5s) at (112.5:1.2) [v] {$r^5\!s$};
      \node (r6s) at (135:1.2) [v] {$r^6\!s$};
      \node (r7s) at (157.5:1.2) [v] {$r^7\!s$};
      \node (r8s) at (180:1.2) [v] {$r^8\!s$};
      \node (r9s) at (-157.5:1.2) [v] {$r^9\!s$};
      \tikzstyle{every node}=[font=\tiny]
      \node (r10s) at (-135:1.2) [v] {$r^{10}\!s$};
      \node (r11s) at (-112.5:1.2) [v] {$r^{11}\!s$};
      \node (r12s) at (-90:1.2) [v] {$r^{12}\!s$};
      \node (r13s) at (-67.5:1.2) [v] {$r^{13}\!s$};
      \node (r14s) at (-45:1.2) [v] {$r^{14}\!s$};
      \node (r15s) at (-22.5:1.2) [v] {$r^{15}\!s$};
      \draw [R] (1) to (r); \draw [R] (r) to (r2); \draw [R] (r2) to (r3);
      \draw [R] (r3) to (r4); \draw [R] (r4) to (r5); \draw [R] (r5) to (r6);
      \draw [R] (r6) to (r7); \draw [R] (r7) to (r8); \draw [R] (r8) to (r9);
      \draw [R] (r9) to (r10); \draw [R] (r10) to (r11);
      \draw [R] (r11) to (r12); \draw [R] (r12) to (r13);
      \draw [R] (r13) to (r14); \draw [R] (r14) to (r15);
      \draw [R] (r15) to (1);
      \draw [r] (s) to (r7s); \draw [r] (r7s) to (r14s);
      \draw [r] (r14s) to (r5s); \draw [r] (r5s) to (r12s);
      \draw [r] (r12s) to (r3s); \draw [r] (r3s) to (r10s);
      \draw [r] (r10s) to (rs); \draw [r] (rs) to (r8s);
      \draw [r] (r8s) to (r15s); \draw [r] (r15s) to (r6s);
      \draw [r] (r6s) to (r13s); \draw [r] (r13s) to (r4s);
      \draw [r] (r4s) to (r11s); \draw [r] (r11s) to (r2s);
      \draw [r] (r2s) to (r9s); \draw [r] (r9s) to (s); 
      \draw [bb] (1) to (s); \draw [bb] (r) to (rs);
      \draw [bb] (r2) to (r2s); \draw [bb] (r3) to (r3s);
      \draw [bb] (r4) to (r4s); \draw [bb] (r5) to (r5s);
      \draw [bb] (r6) to (r6s); \draw [bb] (r7) to (r7s);
      \draw [bb] (r8) to (r8s); \draw [bb] (r9) to (r9s);
      \draw [bb] (r10) to (r10s); \draw [bb] (r11) to (r11s);
      \draw [bb] (r12) to (r12s); \draw [bb] (r13) to (r13s);
      \draw [bb] (r14) to (r14s); \draw [bb] (r15) to (r15s);
      \node at (0,0) {\scriptsize $\SD_{16}$};
    \end{scope}
    \begin{scope}[shift={(5,0)}]
      \tikzstyle{r-bend} = [draw, very thick, eRed,-stealth,bend right=8]
      \node (1) at (0:2) [v] {$1$};
      \node (r) at (22.5:2) [v] {$r$};
      \node (r2) at (45:2) [v] {$r^2$};
      \node (r3) at (67.5:2) [v] {$r^3$};
      \node (r4) at (90:2) [v] {$r^4$};
      \node (r5) at (112.5:2) [v] {$r^5$};
      \node (r6) at (135:2) [v] {$r^6$};
      \node (r7) at (157.5:2) [v] {$r^7$};
      \node (r8) at (180:2) [v] {$r^8$};
      \node (r9) at (-157.5:2) [v] {$r^9$};
      \node (r10) at (-135:2) [v] {$r^{10}$};
      \node (r11) at (-112.5:2) [v] {$r^{11}$};
      \node (r12) at (-90:2) [v] {$r^{12}$};
      \node (r13) at (-67.5:2) [v] {$r^{13}$};
      \node (r14) at (-45:2) [v] {$r^{14}$};
      \node (r15) at (-22.5:2) [v] {$r^{15}$};
      \node (s) at (0:1.2) [v] {$s$};
      \node (rs) at (22.5:1.2) [v] {$r\!s$};
      \node (r2s) at (45:1.2) [v] {$r^2\!s$};
      \node (r3s) at (67.5:1.2) [v] {$r^3\!s$};
      \node (r4s) at (90:1.2) [v] {$r^4\!s$};
      \node (r5s) at (112.5:1.2) [v] {$r^5\!s$};
      \node (r6s) at (135:1.2) [v] {$r^6\!s$};
      \node (r7s) at (157.5:1.2) [v] {$r^7\!s$};
      \node (r8s) at (180:1.2) [v] {$r^8\!s$};
      \node (r9s) at (-157.5:1.2) [v] {$r^9\!s$};
      \tikzstyle{every node}=[font=\tiny]
      \node (r10s) at (-135:1.2) [v] {$r^{10}\!s$};
      \node (r11s) at (-112.5:1.2) [v] {$r^{11}\!s$};
      \node (r12s) at (-90:1.2) [v] {$r^{12}\!s$};
      \node (r13s) at (-67.5:1.2) [v] {$r^{13}\!s$};
      \node (r14s) at (-45:1.2) [v] {$r^{14}\!s$};
      \node (r15s) at (-22.5:1.2) [v] {$r^{15}\!s$};

      \draw [r-bend] (1) to (r); \draw [r-bend] (r) to (r2);
      \draw [r-bend] (r2) to (r3); \draw [r-bend] (r3) to (r4);
      \draw [r-bend] (r4) to (r5); \draw [r-bend] (r5) to (r6);
      \draw [r-bend] (r6) to (r7); \draw [r-bend] (r7) to (r8);
      \draw [r-bend] (r8) to (r9); \draw [r-bend] (r9) to (r10);
      \draw [r-bend] (r10) to (r11); \draw [r-bend] (r11) to (r12);
      \draw [r-bend] (r12) to (r13); \draw [r-bend] (r13) to (r14);
      \draw [r-bend] (r14) to (r15); \draw [r-bend] (r15) to (1);
      \draw [r] (s) to (r9s); \draw [r] (r9s) to (r2s);
      \draw [r] (r2s) to (r11s); \draw [r] (r11s) to (r4s);
      \draw [r] (r4s) to (r13s); \draw [r] (r13s) to (r6s);
      \draw [r] (r6s) to (r15s); \draw [r] (r15s) to (r8s);
      \draw [r] (r8s) to (rs); \draw [r] (rs) to (r10s);
      \draw [r] (r10s) to (r3s); \draw [r] (r3s) to (r12s);
      \draw [r] (r12s) to (r5s); \draw [r] (r5s) to (r14s);
      \draw [r] (r14s) to (r7s); \draw [r] (r7s) to (s); 
      \draw [bb] (1) to (s); \draw [bb] (r) to (rs);
      \draw [bb] (r2) to (r2s); \draw [bb] (r3) to (r3s);
      \draw [bb] (r4) to (r4s); \draw [bb] (r5) to (r5s);
      \draw [bb] (r6) to (r6s); \draw [bb] (r7) to (r7s);
      \draw [bb] (r8) to (r8s); \draw [bb] (r9) to (r9s);
      \draw [bb] (r10) to (r10s); \draw [bb] (r11) to (r11s);
      \draw [bb] (r12) to (r12s); \draw [bb] (r13) to (r13s);
      \draw [bb] (r14) to (r14s); \draw [bb] (r15) to (r15s);
      \node at (0,0) {\scriptsize $\SA_{16}$};
    \end{scope}
  \end{tikzpicture}
  \caption{The four semidirect products of $C_{16}$ with $C_2$.}\label{fig:order32}
\end{figure}

It turns out that for each $n\geq 4$, the generalized quaternion group $Q_{2^n}$ is one of just six groups of order $2^n$ that has an index-$2$ cyclic subgroup. Other than $Q_{2^n}$ and the cyclic group $C_{2^n}$, the other four are semidirect products of $\<r\>\cong C_{2^{n-1}}$ with $C_2$. Each ones\ has a presentation of hte form
\[
\big\<r,s\mid r^{2^{n-1}}=s^2=1,\,srs=r^k\big\>
\]
where $k$ is one of the four solutions to $x^2\equiv 1\pmod{2^n}$: $k=\pm 1$ and $k=2^{n-1}\pm 1$. Note that $k=1$ gives an abelian group $C_{2^{n-1}}\times C_2$, and $k=-1$ defines a dihedral group $D_{2^{n-1}}$. The other two groups are lesser known: $k=2^{n-1}-1$ defines the \emph{semidihedral group} $\SD_{2^{n-1}}$, and $k=2^{n-1}-1$ defined a group that does not even have a standard name. We will call it the \emph{semiabelian group} $\SA_{2^{n-1}}$, which was a name proposed by a student in the author's undergraduate algebra class in 2021. See \cite{macauley2024dihedralizing} for more reasons why this name works well, and more information about both of these groups, such as how to represent them with $2\times 2$ matrices. The Cayley graphs of these last two groups are particularly pleasing, and are shown in Figure~\ref{fig:order32}.

% Chapter 4: proof of only 6 such subgroups
%https://web.mat.bham.ac.uk/D.A.Craven/docs/lectures/pgroups.pdf

%% Semiabelian group for $n=4$ and $n=5%. 

%%%%%%%%%%%%%%%%%%%%%%%%%%%%%%%%%%%%%%%%%%%%%%
\section*{Non-Cayley Graphs and Tables}

Every Cayley graph of a group must have a few basic properties. For one, every node must have one outgoing arrow of each type (color). Additionally, graph also must be \emph{vertex transitive}---every node must look like every other. This forces Cayley graphs to have a structural ``regularity.'' A natural question to ask is whether every Cayley graph with these basic properties must describe a group. This question is vague enough that the answer is ambiguous, but we will explore it now.

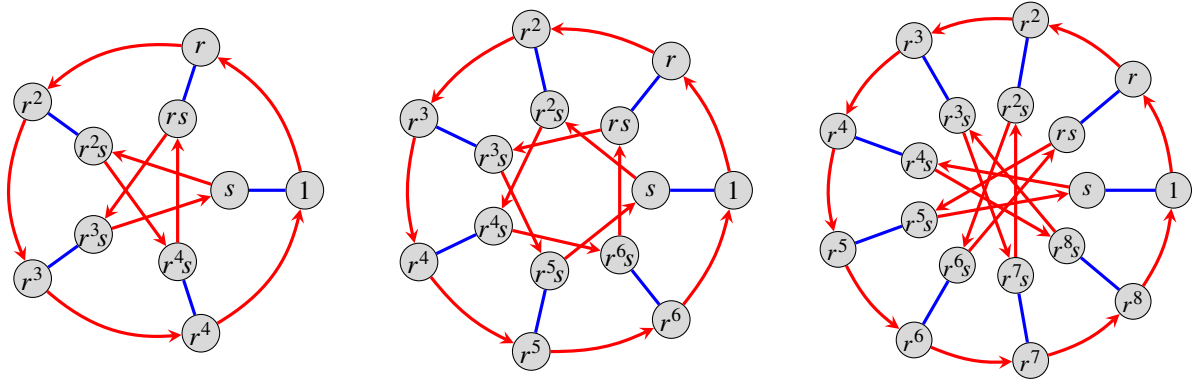
\begin{figure}[!ht] \centering
  \begin{tikzpicture}[scale=1,auto]
  \tikzstyle{v} = [circle, draw, fill=lightgray,inner sep=0pt, 
  minimum size=5mm]
    \tikzstyle{every node}=[font=\small]
    \begin{scope}[shift={(0,0)}] 
    \tikzstyle{r-out} = [draw, very thick, eRed,-stealth,bend right=23]
      \node (e) at (0:2) [v] {$1$};
      \node (r) at (72:2) [v] {$r$};
      \node (r2) at (144:2) [v] {$r^2$};
      \node (r3) at (216:2) [v] {$r^3$};
      \node (r4) at (288:2) [v] {$r^4$};
      \node (f) at (0:1) [v] {$s$};
      \node (rf) at (72:1) [v] {$rs$};
      \node (r2f) at (144:1) [v] {$r^2\!s$};
      \node (r3f) at (216:1) [v] {$r^3\!s$};
      \node (r4f) at (288:1) [v] {$r^4\!s$};
      \draw [r-out] (e) to (r);
      \draw [r-out] (r) to (r2);
      \draw [r-out] (r2) to (r3);
      \draw [r-out] (r3) to (r4);
      \draw [r-out] (r4) to (e);
      \draw [r] (f) to (r2f);
      \draw [r] (r2f) to (r4f);
      \draw [r] (r4f) to (rf);
      \draw [r] (rf) to (r3f);
      \draw [r] (r3f) to (f);
      \draw [bb] (e) to (f);
      \draw [bb] (r) to (rf);
      \draw [bb] (r2) to (r2f);
      \draw [bb] (r3) to (r3f);
      \draw [bb] (r4) to (r4f);
    \end{scope}
      \begin{scope}[shift={(5.5,0)},scale=1.1] 
      \tikzstyle{r-out} = [draw, very thick, eRed,-stealth,bend right=13]
        \node (1) at (0:2) [v] {$1$};
        \node (r) at (51.43:2) [v] {$r$};
        \node (r2) at (102.86:2) [v] {$r^2$};
        \node (r3) at (154.28:2) [v] {$r^3$};
        \node (r4) at (205.71:2) [v] {$r^4$};
        \node (r5) at (257.14:2) [v] {$r^5$};
        \node (r6) at (308.57:2) [v] {$r^6$};
        \node (s) at (0:1) [v] {$s$};
        \node (rs) at (51.43:1) [v] {$rs$};
        \node (r2s) at (102.86:1) [v] {$r^2\!s$};
        \node (r3s) at (154.28:1) [v] {$r^3\!s$};
        \node (r4s) at (205.71:1) [v] {$r^4\!s$};
        \node (r5s) at (257.14:1) [v] {$r^5\!s$};
        \node (r6s) at (308.57:1) [v] {$r^6\!s$};
        \draw [r-out] (1) to (r);
        \draw [r-out] (r) to (r2);
        \draw [r-out] (r2) to (r3);
        \draw [r-out] (r3) to (r4);
        \draw [r-out] (r4) to (r5);
        \draw [r-out] (r5) to (r6);
        \draw [r-out] (r6) to (1);
        \draw [r] (s) to (r2s);
        \draw [r] (r2s) to (r4s);
        \draw [r] (r4s) to (r6s);
        \draw [r] (r6s) to (rs);
        \draw [r] (rs) to (r3s);
        \draw [r] (r3s) to (r5s);
        \draw [r] (r5s) to (s);
        \draw [bb] (1) to (s);
        \draw [bb] (r) to (rs);
        \draw [bb] (r2) to (r2s);
        \draw [bb] (r3) to (r3s);
        \draw [bb] (r4) to (r4s);
        \draw [bb] (r5) to (r5s);
        \draw [bb] (r6) to (r6s);
      \end{scope}
    \begin{scope}[shift={(11.25,0)},scale=1.15] 
      \tikzstyle{v} = [circle, draw, fill=lightgray,inner sep=0pt, 
  minimum size=4.75mm]
     \tikzstyle{every node}=[font=\footnotesize]
      \tikzstyle{r-out} = [draw, very thick, eRed,-stealth,bend right=13]
        \node (1) at (0:2) [v] {$1$};
        \node (r) at (40:2) [v] {$r$};
        \node (r2) at (80:2) [v] {$r^2$};
        \node (r3) at (120:2) [v] {$r^3$};
        \node (r4) at (160:2) [v] {$r^4$};
        \node (r5) at (200:2) [v] {$r^5$};
        \node (r6) at (240:2) [v] {$r^6$};
        \node (r7) at (280:2) [v] {$r^7$};
        \node (r8) at (320:2) [v] {$r^8$};
        \node (s) at (0:1) [v] {$s$};
        \node (rs) at (40:1) [v] {$rs$};
        \node (r2s) at (80:1) [v] {$r^2\!s$};
        \node (r3s) at (120:1) [v] {$r^3\!s$};
        \node (r4s) at (160:1) [v] {$r^4\!s$};
        \node (r5s) at (200:1) [v] {$r^5\!s$};
        \node (r6s) at (240:1) [v] {$r^6\!s$};
        \node (r7s) at (280:1) [v] {$r^7\!s$};
        \node (r8s) at (320:1) [v] {$r^8\!s$};
        \draw [r-out] (1) to (r);
        \draw [r-out] (r) to (r2);
        \draw [r-out] (r2) to (r3);
        \draw [r-out] (r3) to (r4);
        \draw [r-out] (r4) to (r5);
        \draw [r-out] (r5) to (r6);
        \draw [r-out] (r6) to (r7);
        \draw [r-out] (r7) to (r8);
        \draw [r-out] (r8) to (1);
        \draw [r] (s) to (r4s);
        \draw [r] (r4s) to (r8s);
        \draw [r] (r8s) to (r3s);
        \draw [r] (r3s) to (r7s);
        \draw [r] (r7s) to (r2s);
        \draw [r] (r2s) to (r6s);
        \draw [r] (r6s) to (rs);
        \draw [r] (rs) to (r5s);
        \draw [r] (r5s) to (s);
        \draw [bb] (1) to (s);
        \draw [bb] (r) to (rs);
        \draw [bb] (r2) to (r2s);
        \draw [bb] (r3) to (r3s);
        \draw [bb] (r4) to (r4s);
        \draw [bb] (r5) to (r5s);
        \draw [bb] (r6) to (r6s);
        \draw [bb] (r7) to (r7s);
        \draw [bb] (r8) to (r8s);
      \end{scope}
  \end{tikzpicture}
  \caption{Left: a graph with structural regularity that looks like a Cayley graph of a group, except isn't one. Middle and right: two related graphs. Do these define a group?}\label{fig:Petersen}
\end{figure}

Consider the three graphs shown in Figure~\ref{fig:Petersen}, each of which appears to describe a group $G=\<r,s\>$. Note that the skeleton of the first one is the famous Petersen graph. We encourage the reader to spend a few minutes investigating this, and to try to figure out if they describe a group, and if so, which one.

%% Petersen graph non-example
%%

Though the graphs in Figure~\ref{fig:Petersen} look like Cayley graphs, issues appear upon closer inspection. For example, starting with the graph on the left, it is easy to check that it takes ten iterations of the red-blue path to return to the identity element. Therefore, if this is a group, then it must be generated by $rs$, which means it is cyclic. However, cyclic groups are abelian, whereas the red-blue and blue-red paths are different.

All is not lost, though, because this Cayley-ish graph still defines a presentation, and thus a group. But which group? It must be generated by $r$ and $s$, and any path that is a loop defines a relation. Though it is not clear \emph{a priori} how many paths (relations) are needed, let's start with the presentation
\[
\big\<r,s\mid r^5=s^2=1,\,r^3s=sr,\,srs=r^2\big\>.
\]
Since $s^{-1}=s$, the relation $r^3s=sr$ leads to $r^3=srs$, but $srs=r^2$, and so $r^2=r^3$, and hence $r=1$. It is not possible to eliminate $s$ in this manner, and so this graph defines the cyclic group $\<s\mid s^2=1\>\cong C_2$. We encourage the reader to explore the other two examples. The middle example should be quite similar to the first, because $7$, like $5$, is prime. However, the last example has shades of both the first two, and of the semidihedral group in Figure~\ref{fig:order32}. Though there are only two groups of order $14$: $C_{14}$ and $D_7$, there are five groups of order $18$: $C_{18}$, $D_9$, $C_3\times D_3$, $C_3\rtimes D_3$, and $C_3\times C_6$. The interested reader may find it helpful to consult an online searchable group database, like the LMFDB~\cite{lmfdb}.

\begin{figure}[!ht] \centering
  \begin{tikzpicture}[scale=.75,box/.style={anchor=south}]
\colorlet{color-e}{tYellow}
\colorlet{color-a}{tRed}
\colorlet{color-b}{tBlue}
\colorlet{color-c}{tGreen}
\colorlet{color-d}{tPurple}
\newcommand*{\n}{6}%
    \begin{scope}[shift={(0,0)}]
      \path[fill=color-e] (-.1,5) rectangle ++(1,1);
      \path[fill=color-a] (-.1,4) rectangle ++(1,1);
      \path[fill=color-b] (-.1,3) rectangle ++(1,1);
      \path[fill=color-c] (-.1,2) rectangle ++(1,1);
      \path[fill=color-d] (-.1,1) rectangle ++(1,1);
      \path[fill=color-e] (1,6.1) rectangle ++(1,1);
      \path[fill=color-a] (2,6.1) rectangle ++(1,1);
      \path[fill=color-b] (3,6.1) rectangle ++(1,1);
      \path[fill=color-c] (4,6.1) rectangle ++(1,1);
      \path[fill=color-d] (5,6.1) rectangle ++(1,1);
      \path[fill=color-e] (1,5) rectangle ++(1,1);
      \path[fill=color-a] (1,4) rectangle ++(1,1);
      \path[fill=color-b] (1,3) rectangle ++(1,1);
      \path[fill=color-c] (1,2) rectangle ++(1,1);
      \path[fill=color-d] (1,1) rectangle ++(1,1);
      \path[fill=color-a] (2,5) rectangle ++(1,1);
      \path[fill=color-c] (2,4) rectangle ++(1,1);
      \path[fill=color-d] (2,3) rectangle ++(1,1);
      \path[fill=color-b] (2,2) rectangle ++(1,1);
      \path[fill=color-e] (2,1) rectangle ++(1,1);
      \path[fill=color-b] (3,5) rectangle ++(1,1);
      \path[fill=color-d] (3,4) rectangle ++(1,1);
      \path[fill=color-a] (3,3) rectangle ++(1,1);
      \path[fill=color-e] (3,2) rectangle ++(1,1);
      \path[fill=color-c] (3,1) rectangle ++(1,1);
      \path[fill=color-c] (4,5) rectangle ++(1,1);
      \path[fill=color-b] (4,4) rectangle ++(1,1);
      \path[fill=color-e] (4,3) rectangle ++(1,1);
      \path[fill=color-d] (4,2) rectangle ++(1,1);
      \path[fill=color-a] (4,1) rectangle ++(1,1);
      \path[fill=color-d] (5,5) rectangle ++(1,1);
      \path[fill=color-e] (5,4) rectangle ++(1,1);
      \path[fill=color-c] (5,3) rectangle ++(1,1);
      \path[fill=color-a] (5,2) rectangle ++(1,1);
      \path[fill=color-b] (5,1) rectangle ++(1,1);
      \foreach \i in {1,...,\n} {
        \draw [very thin] (\i,1) -- (\i,\n);  % vertical lines (columns)
        \draw [very thin] (1,\i) -- (\n,\i);  % horizontal lines (rows)
        \draw [very thin] (\i,\n+.1) -- (\i,\n+1.1);  % short vertical lines 
        \draw [very thin] (-.1,\i) -- (.9,\i);  % short horizontal lines
      }
      \draw [very thin] (-.1,1) -- (-.1,\n);
      \draw [very thin] (.9,1) -- (.9,\n);
      \draw [very thin] (1,\n+.1) -- (\n,\n+.1);
      \draw [very thin] (1,\n+1.1) -- (\n,\n+1.1);
      %%
      %% Column headers
      \node [box] at (0.4,5.15) {$e$};
      \node [box] at (0.4,4.15) {$a$};
      \node [box] at (0.4,3.15) {$b$};
      \node [box] at (0.4,2.15) {$c$};
      \node [box] at (0.4,1.15) {$d$};
      %%
      %% Row headers 
      \node [box] at (1.5,6.25) {$e$};
      \node [box] at (2.5,6.25) {$a$};
      \node [box] at (3.5,6.25) {$b$};
      \node [box] at (4.5,6.25) {$c$};
      \node [box] at (5.5,6.25) {$d$};
      %%
      % 1st Column: e
      \node [box] at (1.5,5.15) {$e$};
      \node [box] at (1.5,4.15) {$a$};
      \node [box] at (1.5,3.15) {$b$};
      \node [box] at (1.5,2.15) {$c$};
      \node [box] at (1.5,1.15) {$d$};
      % 2nd Column: a
      \node [box] at (2.5,5.15) {$a$};
      \node [box] at (2.5,4.15) {$c$};
      \node [box] at (2.5,3.15) {$d$};
      \node [box] at (2.5,2.15) {$b$};
      \node [box] at (2.5,1.15) {$e$};
      % 3rd Column: b
      \node [box] at (3.5,5.15) {$b$};
      \node [box] at (3.5,4.15) {$d$};
      \node [box] at (3.5,3.15) {$a$};
      \node [box] at (3.5,2.15) {$e$};
      \node [box] at (3.5,1.15) {$c$};
      % 4th Column: c
      \node [box] at (4.5,5.15) {$c$};
      \node [box] at (4.5,4.15) {$b$};
      \node [box] at (4.5,3.15) {$e$};
      \node [box] at (4.5,2.15) {$d$};
      \node [box] at (4.5,1.15) {$a$};
      % 5th Column: d
      \node [box] at (5.5,5.15) {$d$};
      \node [box] at (5.5,4.15) {$e$};
      \node [box] at (5.5,3.15) {$c$};
      \node [box] at (5.5,2.15) {$a$};
      \node [box] at (5.5,1.15) {$b$};
    \end{scope}
    \begin{scope}[shift={(10,0)}]
      \path[fill=color-e] (-.1,5) rectangle ++(1,1);
      \path[fill=color-a] (-.1,4) rectangle ++(1,1);
      \path[fill=color-b] (-.1,3) rectangle ++(1,1);
      \path[fill=color-c] (-.1,2) rectangle ++(1,1);
      \path[fill=color-d] (-.1,1) rectangle ++(1,1);
      \path[fill=color-e] (1,6.1) rectangle ++(1,1);
      \path[fill=color-a] (2,6.1) rectangle ++(1,1);
      \path[fill=color-b] (3,6.1) rectangle ++(1,1);
      \path[fill=color-c] (4,6.1) rectangle ++(1,1);
      \path[fill=color-d] (5,6.1) rectangle ++(1,1);
      \path[fill=color-e] (1,5) rectangle ++(1,1);
      \path[fill=color-a] (1,4) rectangle ++(1,1);
      \path[fill=color-b] (1,3) rectangle ++(1,1);
      \path[fill=color-c] (1,2) rectangle ++(1,1);
      \path[fill=color-d] (1,1) rectangle ++(1,1);
      \path[fill=color-a] (2,5) rectangle ++(1,1);
      \path[fill=color-e] (2,4) rectangle ++(1,1);
      \path[fill=color-d] (2,3) rectangle ++(1,1);
      \path[fill=color-b] (2,2) rectangle ++(1,1);
      \path[fill=color-c] (2,1) rectangle ++(1,1);
      \path[fill=color-b] (3,5) rectangle ++(1,1);
      \path[fill=color-c] (3,4) rectangle ++(1,1);
      \path[fill=color-e] (3,3) rectangle ++(1,1);
      \path[fill=color-d] (3,2) rectangle ++(1,1);
      \path[fill=color-a] (3,1) rectangle ++(1,1);
      \path[fill=color-c] (4,5) rectangle ++(1,1);
      \path[fill=color-d] (4,4) rectangle ++(1,1);
      \path[fill=color-a] (4,3) rectangle ++(1,1);
      \path[fill=color-e] (4,2) rectangle ++(1,1);
      \path[fill=color-b] (4,1) rectangle ++(1,1);
      \path[fill=color-d] (5,5) rectangle ++(1,1);
      \path[fill=color-b] (5,4) rectangle ++(1,1);
      \path[fill=color-c] (5,3) rectangle ++(1,1);
      \path[fill=color-a] (5,2) rectangle ++(1,1);
      \path[fill=color-e] (5,1) rectangle ++(1,1);
      \foreach \i in {1,...,\n} {
        \draw [very thin] (\i,1) -- (\i,\n);  % vertical lines (columns)
        \draw [very thin] (1,\i) -- (\n,\i);  % horizontal lines (rows)
        \draw [very thin] (\i,\n+.1) -- (\i,\n+1.1);  % short vertical lines 
        \draw [very thin] (-.1,\i) -- (.9,\i);  % short horizontal lines
      }
      \draw [very thin] (-.1,1) -- (-.1,\n);
      \draw [very thin] (.9,1) -- (.9,\n);
      \draw [very thin] (1,\n+.1) -- (\n,\n+.1);
      \draw [very thin] (1,\n+1.1) -- (\n,\n+1.1);
      %%
      %% Column headers
      \node [box] at (0.4,5.15) {$e$};
      \node [box] at (0.4,4.15) {$a$};
      \node [box] at (0.4,3.15) {$b$};
      \node [box] at (0.4,2.15) {$c$};
      \node [box] at (0.4,1.15) {$d$};
      %%
      %% Row headers 
      \node [box] at (1.5,6.25) {$e$};
      \node [box] at (2.5,6.25) {$a$};
      \node [box] at (3.5,6.25) {$b$};
      \node [box] at (4.5,6.25) {$c$};
      \node [box] at (5.5,6.25) {$d$};
      %%
      % 1st Column: e
      \node [box] at (1.5,5.15) {$e$};
      \node [box] at (1.5,4.15) {$a$};
      \node [box] at (1.5,3.15) {$b$};
      \node [box] at (1.5,2.15) {$c$};
      \node [box] at (1.5,1.15) {$d$};
      % 2nd Column: a
      \node [box] at (2.5,5.15) {$a$};
      \node [box] at (2.5,4.15) {$e$};
      \node [box] at (2.5,3.15) {$d$};
      \node [box] at (2.5,2.15) {$b$};
      \node [box] at (2.5,1.15) {$c$};
      % 3rd Column: b
      \node [box] at (3.5,5.15) {$b$};
      \node [box] at (3.5,4.15) {$c$};
      \node [box] at (3.5,3.15) {$e$};
      \node [box] at (3.5,2.15) {$d$};
      \node [box] at (3.5,1.15) {$a$};
      % 4th Column: c
      \node [box] at (4.5,5.15) {$c$};
      \node [box] at (4.5,4.15) {$d$};
      \node [box] at (4.5,3.15) {$a$};
      \node [box] at (4.5,2.15) {$e$};
      \node [box] at (4.5,1.15) {$b$};
      % 5th Column: d
      \node [box] at (5.5,5.15) {$d$};
      \node [box] at (5.5,4.15) {$b$};
      \node [box] at (5.5,3.15) {$c$};
      \node [box] at (5.5,2.15) {$a$};
      \node [box] at (5.5,1.15) {$e$};
    \end{scope}
  \end{tikzpicture}
  \caption{Two Latin squares. Do these represent groups?}\label{fig:latin-squares}
\end{figure}
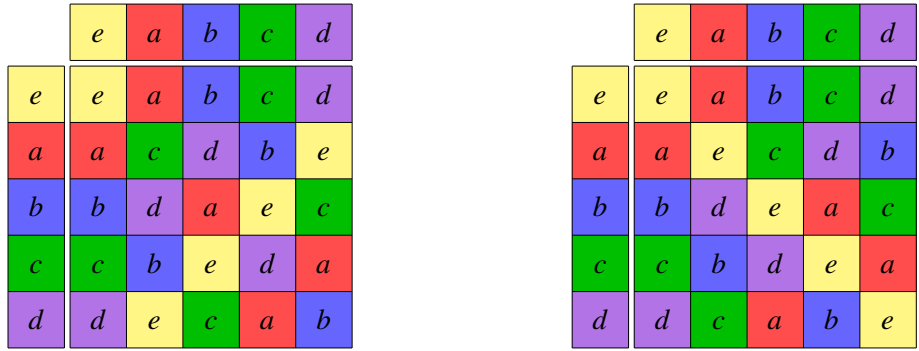

Just like how we asked if every graph that ``looks like'' a Cayley graph is one, we can ask if every table that looks like a Cayley table defines a group. It is easy to translate the requirement that every group has an identity element into a basic property of a table---the (without loss of generality) first row and column must agree with the headers. Beyond this, every row and column must contain each element exactly once. Such a table is called a \emph{Latin square}. Figure~\ref{fig:latin-squares} shows two Latin squares of order $5$. We encourage the reader to ponder whether either of these form a group, before preceding on, where  the answer is given. 

%% Latin square non-example
%%

The first Latin square in Figure~\ref{fig:latin-squares} describes the (additive) cyclic group $\Z_5=\{0,1,2,3,4\}$, with $e=0$, $a=1$, $b=3$, $c=2$, and $d=4$. The one on the right cannot be a group, because $g^2=e$ would imply that every non-identity element has order $2$, which is impossible in a group of odd order, due to Lagrange's theorem. Upon closer inspection, associativity fails for this table, because $(a*b)*d=c*d=a$, but $a*(b*d)=a*c=d$.

%%%%%%%%%%%%%%%%%%%%%%%%%%%%%%%%%%%%%%%%%%%%%%
\section*{Exploratory Examples}

Now that we have seen several examples and non-examples of groups, we will finish with some diagrams that may or may not be Cayley graphs. As before, all of these define a presentation, and so they define groups. However, if they are not Cayley graphs, then the groups they define are smaller than the sizes of the vertex sets. We will not actually give the answers to any of these, because we want to leave a few tantalizing puzzles for the curious reader. Anyone who wants to check their answers can use the open-source GAP (Groups, Algorithms and Programming) software, which can take in a presentation defined by such a diagram, and determine which group it is. Note that this can be cross-checked with an online group database such as the LMFDB. Finally, we want to encourage interested readers to create and share their own puzzles like this. Even undergraduate algebra students can create these types of diagrams, and doing so could be a fun project.

\begin{figure}[!ht] \centering
  \begin{tikzpicture}[scale=1.6,auto]
  \tikzstyle{v} = [circle, draw, fill=lightgray,inner sep=0pt, 
    minimum size=4mm]
    \tikzstyle{every node}=[font=\footnotesize]
          \begin{scope}[shift={(0,0)}]
    \tikzstyle{rr-bend} = [draw, very thick, eRed,bend right=13]
    \tikzstyle{bb-bend} = [draw, very thick, eBlue,bend right=13]
        \node (s) at (0:1.1) [v] {};
        \node (rs) at (45:1.15) [v] {};
        \node (r2s) at (90:1.15) [v] {};
        \node (r3s) at (135:1.15) [v] {};
        \node (r4s) at (180:1.15) [v] {};
        \node (r5s) at (225:1.15) [v] {};
        \node (r6s) at (270:1.15) [v] {};
        \node (r7s) at (315:1.15) [v] {};
        \node (1) at (0:2) [v] {$I$};
        \node (r) at (45:2) [v] {};
        \node (r2) at (90:2) [v] {};
        \node (r3) at (135:2) [v] {};
        \node (r4) at (180:2) [v] {};
        \node (r5) at (225:2) [v] {};
        \node (r6) at (270:2) [v] {};
        \node (r7) at (315:2) [v] {};
        \draw [bb-bend] (1) to (r);
        \draw [rr-bend] (r) to (r2);
        \draw [bb-bend] (r2) to (r3);
        \draw [rr-bend] (r3) to (r4);
        \draw [bb-bend] (r4) to (r5);
        \draw [rr-bend] (r5) to (r6);
        \draw [bb-bend] (r6) to (r7);
        \draw [rr-bend] (r7) to (1);
        \draw [bb] (s) to (r3s);
        \draw [rr] (r3s) to (r6s);
        \draw [bb] (r6s) to (rs);
        \draw [rr] (rs) to (r4s);
        \draw [bb] (r4s) to (r7s);
        \draw [rr] (r7s) to (r2s);
        \draw [bb] (r2s) to (r5s);
        \draw [rr] (r5s) to (s);
        \draw [gg] (1) to (s); \draw [gg] (r) to (rs);
        \draw [gg] (r2) to (r2s); \draw [gg] (r3) to (r3s);
        \draw [gg] (r4) to (r4s); \draw [gg] (r5) to (r5s);
        \draw [gg] (r6) to (r6s); \draw [gg] (r7) to (r7s);
      \end{scope}
    \begin{scope}[shift={(5,0)}]
    \tikzstyle{rr-bend} = [draw, very thick, eRed,bend left=8]
    \tikzstyle{bb-bend} = [draw, very thick, eBlue,bend left=8]
      \node (s) at (90:1.2) [v] {};
      \node (rs) at (67.5:1.2) [v] {};
      \node (r2s) at (45:1.2) [v] {};
      \node (r3s) at (22.5:1.2) [v] {};
      \node (r4s) at (0:1.2) [v] {};
      \node (r5s) at (-22.5:1.2) [v] {};
      \node (r6s) at (-45:1.2) [v] {};
      \node (r7s) at (-67.5:1.2) [v] {};
      \node (r8s) at (-90:1.2) [v] {};
      \node (r9s) at (-112.5:1.2) [v] {};
      \node (r10s) at (-135:1.2) [v] {};
      \node (r11s) at (-157.5:1.2) [v] {};
      \node (r12s) at (180:1.2) [v] {};
      \node (r13s) at (157.5:1.2) [v] {};
      \node (r14s) at (135:1.2) [v] {};
      \node (r15s) at (112.5:1.2) [v] {};
      \node (1) at (90:2) [v] {};
      \node (r) at (67.5:2) [v] {};
      \node (r2) at (45:2) [v] {};
      \node (r3) at (22.5:2) [v] {};
      \node (r4) at (0:2) [v] {$I$};
      \node (r5) at (-22.5:2) [v] {};
      \node (r6) at (-45:2) [v] {};
      \node (r7) at (-67.5:2) [v] {};
      \node (r8) at (-90:2) [v] {};
      \node (r9) at (-112.5:2) [v] {};
      \node (r10) at (-135:2) [v] {};
      \node (r11) at (-157.5:2) [v] {};
      \node (r12) at (180:2) [v] {};
      \node (r13) at (157.5:2) [v] {};
      \node (r14) at (135:2) [v] {};
      \node (r15) at (112.5:2) [v] {};
      \draw [rr-bend] (1) to (r); \draw [bb-bend] (r) to (r2);
      \draw [rr-bend] (r2) to (r3); \draw [bb-bend] (r3) to (r4);
      \draw [rr-bend] (r4) to (r5); \draw [bb-bend] (r5) to (r6);
      \draw [rr-bend] (r6) to (r7); \draw [bb-bend] (r7) to (r8);
      \draw [rr-bend] (r8) to (r9); \draw [bb-bend] (r9) to (r10);
      \draw [rr-bend] (r10) to (r11); \draw [bb-bend] (r11) to (r12);
      \draw [rr-bend] (r12) to (r13); \draw [bb-bend] (r13) to (r14);
      \draw [rr-bend] (r14) to (r15); \draw [bb-bend] (r15) to (1);
      \draw [bb] (s) to (r9s); \draw [rr] (r9s) to (r2s);
      \draw [bb] (r2s) to (r11s); \draw [rr] (r11s) to (r4s);
      \draw [bb] (r4s) to (r13s); \draw [rr] (r13s) to (r6s);
      \draw [bb] (r6s) to (r15s); \draw [rr] (r15s) to (r8s);
      \draw [bb] (r8s) to (rs); \draw [rr] (rs) to (r10s);
      \draw [bb] (r10s) to (r3s); \draw [rr] (r3s) to (r12s);
      \draw [bb] (r12s) to (r5s); \draw [rr] (r5s) to (r14s);
      \draw [bb] (r14s) to (r7s); \draw [rr] (r7s) to (s); 
      \draw [gg] (1) to (s); \draw [gg] (r) to (rs);
      \draw [gg] (r2) to (r2s); \draw [gg] (r3) to (r3s);
      \draw [gg] (r4) to (r4s); \draw [gg] (r5) to (r5s);
      \draw [gg] (r6) to (r6s); \draw [gg] (r7) to (r7s);
      \draw [gg] (r8) to (r8s); \draw [gg] (r9) to (r9s);
      \draw [gg] (r10) to (r10s); \draw [gg] (r11) to (r11s);
      \draw [gg] (r12) to (r12s); \draw [gg] (r13) to (r13s);
      \draw [gg] (r14) to (r14s); \draw [gg] (r15) to (r15s);
    \end{scope}
  \end{tikzpicture}
  \caption{Cayley graphs that are not quite the same as our standard ones for
    the diquaternion groups, $\DQ_8$ and $\DQ_{16}$.}\label{fig:almost-DQ}
\end{figure}
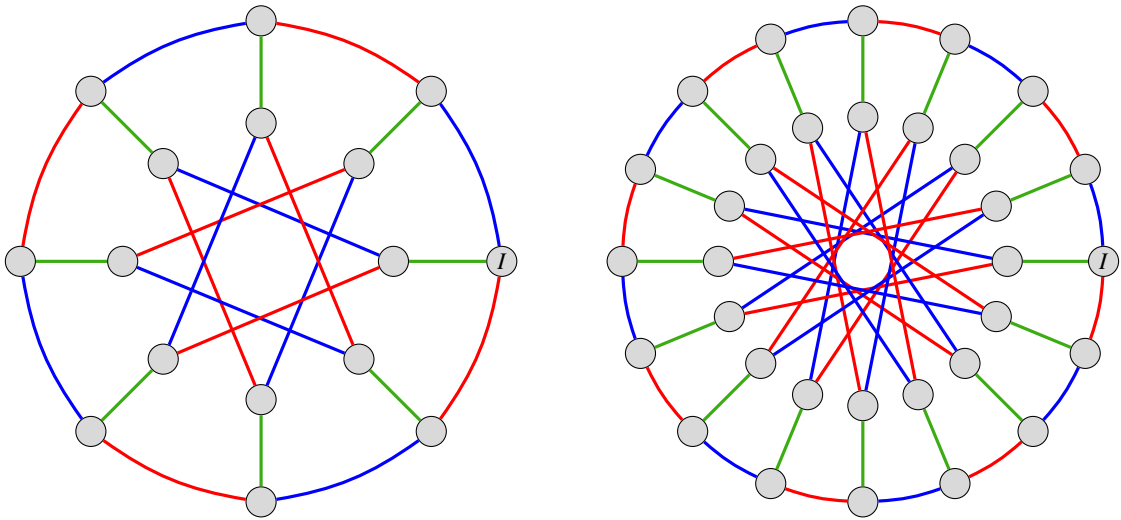

Our first example consists of a modified version of the Cayley graphs of the diquaternion groups, shown in Figure~\ref{fig:almost-DQ}. Note that in $\DQ_8$, the product of any two distinct generators has order $4$, because each of the paths: red-blue-red-blue, red-green-red-green, and blue-green-blue-green, ends up at the antipodal node, $-1$. In contrast, in the graph on the left in Figure~\ref{fig:almost-DQ}, it takes four iterations of each ``\emph{color1-color2} path'' to reach the antipodal node. On the right, each of these takes eight iterations. 

In other words, if the graphs shown in Figure~\ref{fig:almost-DQ} are Cayley graphs, and the generators are denoted $X$, $Y$, and $Z$, then $|XY|=|XZ|=|YZ|=8$ in the group on the left, and $|XY|=|XZ|=|YZ|=16$ in the group on the right. This means that both of these groups have an index-$2$ cyclic subgroup. We have already mentioned that there are only six such groups of order $2^n$ for each $n\geq 4$, four of which are nonabelian. In other words, if the graph on the left is a Cayley graph, then it must be either $D_8$, $\SD_8$, $\SA_8$, or $Q_{16}$. Analogously, the one on the right must define either $D_{16}$, $\SD_{16}$, $\SA_{16}$, or $Q_{32}$, or not be a Cayley graph at all.

\begin{figure}[!ht] \centering
  \begin{tikzpicture}[scale=1.25,auto]
  \tikzstyle{v} = [circle, draw, fill=lightgray,inner sep=0pt, minimum size=3mm]
\tikzstyle{R} = [draw, very thick, eRed,-stealth,bend right=10]
\tikzstyle{R-out} = [draw, very thick, eRed,-stealth,bend right=15]
\tikzstyle{R-in} = [draw, very thick, eRed,-stealth,bend left=12]
\tikzstyle{R2-in} = [draw, very thick, eRed,-stealth,bend right=12]
\tikzstyle{B} = [draw, very thick, eBlue,-stealth,bend right=35]
    %% not a group
    \begin{scope}[shift={(6,0)}]
      \node (s) at (0:.8) [v] {};
      \node (rs) at (45:.8) [v] {};
      \node (r2s) at (90:.8) [v] {};
      \node (r3s) at (135:.8) [v] {};
      \node (r4s) at (180:.8) [v] {};
      \node (r5s) at (225:.8) [v] {};
      \node (r6s) at (270:.8) [v] {};
      \node (r7s) at (315:.8) [v] {};
      \node (1) at (0:2) [v] {};
      \node (r) at (45:2) [v] {};
      \node (r2) at (90:2) [v] {};
      \node (r3) at (135:2) [v] {};
      \node (r4) at (180:2) [v] {};
      \node (r5) at (225:2) [v] {};
      \node (r6) at (270:2) [v] {};
      \node (r7) at (315:2) [v] {};
      \draw [R-out] (1) to (r);
      \draw [R-out] (r) to (r2);
      \draw [R-out] (r2) to (r3);
      \draw [R-out] (r3) to (r4);
      \draw [R-out] (r4) to (r5);
      \draw [R-out] (r5) to (r6);
      \draw [R-out] (r6) to (r7);
      \draw [R-out] (r7) to (1);
      \draw [B] (1) to (rs); \draw [B] (r) to (r2s);
      \draw [B] (r2) to (r3s); \draw [B] (r3) to (r4s);
      \draw [B] (r4) to (r5s); \draw [B] (r5) to (r6s);
      \draw [B] (r6) to (r7s); \draw [B] (r7) to (s);
      \draw [B] (s) to (r); \draw [B] (rs) to (r2);
      \draw [B] (r2s) to (r3); \draw [B] (r3s) to (r4);
      \draw [B] (r4s) to (r5); \draw [B] (r5s) to (r6);
      \draw [B] (r6s) to (r7); \draw [B] (r7s) to (1);
      \draw [R-in] (s) to (r7s);
      \draw [R-in] (r7s) to (r6s);
      \draw [R-in] (r6s) to (r5s);
      \draw [R-in] (r5s) to (r4s);
      \draw [R-in] (r4s) to (r3s);
      \draw [R-in] (r3s) to (r2s);
      \draw [R-in] (r2s) to (rs);
      \draw [R-in] (rs) to (s);
   \end{scope}
    \begin{scope}[shift={(0,0)}]
          \node (s) at (0:.8) [v] {};
      \node (rs) at (45:.8) [v] {};
      \node (r2s) at (90:.8) [v] {};
      \node (r3s) at (135:.8) [v] {};
      \node (r4s) at (180:.8) [v] {};
      \node (r5s) at (225:.8) [v] {};
      \node (r6s) at (270:.8) [v] {};
      \node (r7s) at (315:.8) [v] {};
      \node (1) at (0:2) [v] {};
      \node (r) at (45:2) [v] {};
      \node (r2) at (90:2) [v] {};
      \node (r3) at (135:2) [v] {};
      \node (r4) at (180:2) [v] {};
      \node (r5) at (225:2) [v] {};
      \node (r6) at (270:2) [v] {};
      \node (r7) at (315:2) [v] {};
      \draw [R-out] (1) to (r);
      \draw [R-out] (r) to (r2);
      \draw [R-out] (r2) to (r3);
      \draw [R-out] (r3) to (r4);
      \draw [R-out] (r4) to (r5);
      \draw [R-out] (r5) to (r6);
      \draw [R-out] (r6) to (r7);
      \draw [R-out] (r7) to (1);
      \draw [B] (1) to (rs); \draw [B] (r) to (r2s);
      \draw [B] (r2) to (r3s); \draw [B] (r3) to (r4s);
      \draw [B] (r4) to (r5s); \draw [B] (r5) to (r6s);
      \draw [B] (r6) to (r7s); \draw [B] (r7) to (s);
      \draw [B] (s) to (r); \draw [B] (rs) to (r2);
      \draw [B] (r2s) to (r3); \draw [B] (r3s) to (r4);
      \draw [B] (r4s) to (r5); \draw [B] (r5s) to (r6);
      \draw [B] (r6s) to (r7); \draw [B] (r7s) to (1);
      \draw [R2-in] (r7s) to (s);
      \draw [R2-in] (r6s) to (r7s);
      \draw [R2-in] (r5s) to (r6s);
      \draw [R2-in] (r4s) to (r5s);
      \draw [R2-in] (r3s) to (r4s);
      \draw [R2-in] (r2s) to (r3s);
      \draw [R2-in] (rs) to (r2s);
      \draw [R2-in] (s) to (rs);
   \end{scope}
\end{tikzpicture}
  \caption{Two other graphs on $16$ nodes. Are these the Cayley graphs of a group?}\label{fig:flower-16}
\end{figure}
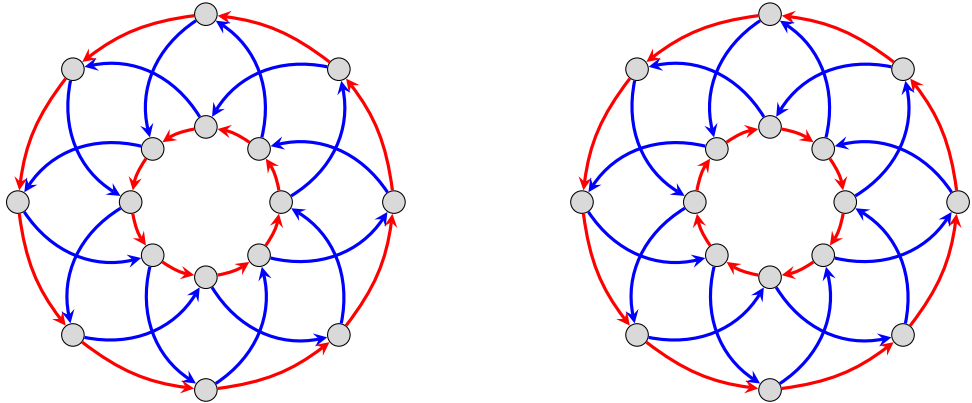

Another family of curious examples is shown in Figure~\ref{fig:flower-16}. Notice how the one on the left, if it is a Cayley graph, must be of an abelian group. Therefore, it must be either $C_{16}$ or $C_8\times C_2$. In contrast, the one on the right, if it is a Cayley graph, must be either $D_8$, $\SD_8$, $\SA_8$, or $Q_{16}$.

%https://math.stackexchange.com/questions/2011666/a-reference-for-a-generalized-quaternion-group-has-a-unique-element-of-order-2

%% Graphs with 16 nodes. 
  
\begin{figure}[!ht]
  \begin{tikzpicture}[scale=1.6]
  \tikzstyle{v} = [circle, draw, fill=lightgray,inner sep=0pt, 
  minimum size=2.5mm]
\tikzstyle{R} = [draw, very thick, eRed,-stealth,bend right=10]
\tikzstyle{B2} = [draw, very thick, eBlue,-stealth,bend right=35]
    \begin{scope}[shift={(0,0)}]
      \node (s) at (0:1) [v] {};
      \node (rs) at (22.5:1) [v] {};
      \node (r2s) at (45:1) [v] {};
      \node (r3s) at (67.5:1) [v] {};
      \node (r4s) at (90:1) [v] {};
      \node (r5s) at (112.5:1) [v] {};
      \node (r6s) at (135:1) [v] {};
      \node (r7s) at (157.5:1) [v] {};
      \node (r8s) at (180:1) [v] {};
      \node (r9s) at (202.5:1) [v] {};
      \node (r10s) at (225:1) [v] {};
      \node (r11s) at (247.5:1) [v] {};
      \node (r12s) at (270:1) [v] {};
      \node (r13s) at (292.5:1) [v] {};
      \node (r14s) at (315:1) [v] {};
      \node (r15s) at (337.5:1) [v] {};
      \node (1) at (0:2) [v] {};
      \node (r) at (22.5:2) [v] {};
      \node (r2) at (45:2) [v] {};
      \node (r3) at (67.5:2) [v] {};
      \node (r4) at (90:2) [v] {};
      \node (r5) at (112.5:2) [v] {};
      \node (r6) at (135:2) [v] {};
      \node (r7) at (157.5:2) [v] {};
      \node (r8) at (180:2) [v] {};
      \node (r9) at (202.5:2) [v] {};
      \node (r10) at (225:2) [v] {};
      \node (r11) at (247.5:2) [v] {};
      \node (r12) at (270:2) [v] {};
      \node (r13) at (292.5:2) [v] {};
      \node (r14) at (315:2) [v] {};
      \node (r15) at (337.5:2) [v] {};
      \draw [r] (s) to (r15s); \draw [r] (r15s) to (r14s);
      \draw [r] (r14s) to (r13s); \draw [r] (r13s) to (r12s);
      \draw [r] (r12s) to (r11s); \draw [r] (r11s) to (r10s);
      \draw [r] (r10s) to (r9s); \draw [r] (r9s) to (r8s);
      \draw [r] (r8s) to (r7s); \draw [r] (r7s) to (r6s);
      \draw [r] (r6s) to (r5s); \draw [r] (r5s) to (r4s);
      \draw [r] (r4s) to (r3s); \draw [r] (r3s) to (r2s);
      \draw [r] (r2s) to (rs); \draw [r] (rs) to (s); 
      \draw [B2] (1) to (r3s); \draw [B2] (r) to (r4s);
      \draw [B2] (r2) to (r5s); \draw [B2] (r3) to (r6s);
      \draw [B2] (r4) to (r7s); \draw [B2] (r5) to (r8s);
      \draw [B2] (r6) to (r9s); \draw [B2] (r7) to (r10s);
      \draw [B2] (r8) to (r11s); \draw [B2] (r9) to (r12s);
      \draw [B2] (r10) to (r13s); \draw [B2] (r11) to (r14s);
      \draw [B2] (r12) to (r15s); \draw [B2] (r13) to (s);
      \draw [B2] (r14) to (rs); \draw [B2] (r15) to (r2s);
      \draw [B2] (s) to (r3); \draw [B2] (rs) to (r4);
      \draw [B2] (r2s) to (r5); \draw [B2] (r3s) to (r6);
      \draw [B2] (r4s) to (r7); \draw [B2] (r5s) to (r8);
      \draw [B2] (r6s) to (r9); \draw [B2] (r7s) to (r10);
      \draw [B2] (r8s) to (r11); \draw [B2] (r9s) to (r12);
      \draw [B2] (r10s) to (r13); \draw [B2] (r11s) to (r14);
      \draw [B2] (r12s) to (r15); \draw [B2] (r13s) to (1);
      \draw [B2] (r14s) to (r); \draw [B2] (r15s) to (r2);
      \draw [R] (1) to (r); \draw [R] (r) to (r2); \draw [R] (r2) to (r3);
      \draw [R] (r3) to (r4); \draw [R] (r4) to (r5); \draw [R] (r5) to (r6);
      \draw [R] (r6) to (r7); \draw [R] (r7) to (r8); \draw [R] (r8) to (r9);
      \draw [R] (r9) to (r10); \draw [R] (r10) to (r11);
      \draw [R] (r11) to (r12); \draw [R] (r12) to (r13);
      \draw [R] (r13) to (r14); \draw [R] (r14) to (r15);
      \draw [R] (r15) to (1);
    \end{scope}
    \begin{scope}[shift={(6,0)}]
      \node (s) at (0:1) [v] {};
      \node (rs) at (22.5:1) [v] {};
      \node (r2s) at (45:1) [v] {};
      \node (r3s) at (67.5:1) [v] {};
      \node (r4s) at (90:1) [v] {};
      \node (r5s) at (112.5:1) [v] {};
      \node (r6s) at (135:1) [v] {};
      \node (r7s) at (157.5:1) [v] {};
      \node (r8s) at (180:1) [v] {};
      \node (r9s) at (202.5:1) [v] {};
      \node (r10s) at (225:1) [v] {};
      \node (r11s) at (247.5:1) [v] {};
      \node (r12s) at (270:1) [v] {};
      \node (r13s) at (292.5:1) [v] {};
      \node (r14s) at (315:1) [v] {};
      \node (r15s) at (337.5:1) [v] {};
      \node (1) at (0:2) [v] {};
      \node (r) at (22.5:2) [v] {};
      \node (r2) at (45:2) [v] {};
      \node (r3) at (67.5:2) [v] {};
      \node (r4) at (90:2) [v] {};
      \node (r5) at (112.5:2) [v] {};
      \node (r6) at (135:2) [v] {};
      \node (r7) at (157.5:2) [v] {};
      \node (r8) at (180:2) [v] {};
      \node (r9) at (202.5:2) [v] {};
      \node (r10) at (225:2) [v] {};
      \node (r11) at (247.5:2) [v] {};
      \node (r12) at (270:2) [v] {};
      \node (r13) at (292.5:2) [v] {};
      \node (r14) at (315:2) [v] {};
      \node (r15) at (337.5:2) [v] {};
      \draw [r] (s) to (r15s); \draw [r] (r15s) to (r14s);
      \draw [r] (r14s) to (r13s); \draw [r] (r13s) to (r12s);
      \draw [r] (r12s) to (r11s); \draw [r] (r11s) to (r10s);
      \draw [r] (r10s) to (r9s); \draw [r] (r9s) to (r8s);
      \draw [r] (r8s) to (r7s); \draw [r] (r7s) to (r6s);
      \draw [r] (r6s) to (r5s); \draw [r] (r5s) to (r4s);
      \draw [r] (r4s) to (r3s); \draw [r] (r3s) to (r2s);
      \draw [r] (r2s) to (rs); \draw [r] (rs) to (s); 
      \draw [B2] (1) to (r2s); \draw [B2] (r) to (r3s);
      \draw [B2] (r2) to (r4s); \draw [B2] (r3) to (r5s);
      \draw [B2] (r4) to (r6s); \draw [B2] (r5) to (r7s);
      \draw [B2] (r6) to (r8s); \draw [B2] (r7) to (r9s);
      \draw [B2] (r8) to (r10s); \draw [B2] (r9) to (r11s);
      \draw [B2] (r10) to (r12s); \draw [B2] (r11) to (r13s);
      \draw [B2] (r12) to (r14s); \draw [B2] (r13) to (r15s);
      \draw [B2] (r14) to (s); \draw [B2] (r15) to (rs);
      \draw [B2] (s) to (r2); \draw [B2] (rs) to (r3);
      \draw [B2] (r2s) to (r4); \draw [B2] (r3s) to (r5);
      \draw [B2] (r4s) to (r6); \draw [B2] (r5s) to (r7);
      \draw [B2] (r6s) to (r8); \draw [B2] (r7s) to (r9);
      \draw [B2] (r8s) to (r10); \draw [B2] (r9s) to (r11);
      \draw [B2] (r10s) to (r12); \draw [B2] (r11s) to (r13);
      \draw [B2] (r12s) to (r14); \draw [B2] (r13s) to (r15);
      \draw [B2] (r14s) to (1); \draw [B2] (r15s) to (r);
      \draw [R] (1) to (r); \draw [R] (r) to (r2); \draw [R] (r2) to (r3);
      \draw [R] (r3) to (r4); \draw [R] (r4) to (r5); \draw [R] (r5) to (r6);
      \draw [R] (r6) to (r7); \draw [R] (r7) to (r8); \draw [R] (r8) to (r9);
      \draw [R] (r9) to (r10); \draw [R] (r10) to (r11);
      \draw [R] (r11) to (r12); \draw [R] (r12) to (r13);
      \draw [R] (r13) to (r14); \draw [R] (r14) to (r15);
      \draw [R] (r15) to (1);
    \end{scope}
  \end{tikzpicture}
  \caption{If either of these describe a group of order $32$, then they must be the Cayley graph of either  $D_{16}$, $\SD_{16}$, $\SA_{16}$, or $Q_{32}$.}\label{fig:flower-32}
\end{figure}
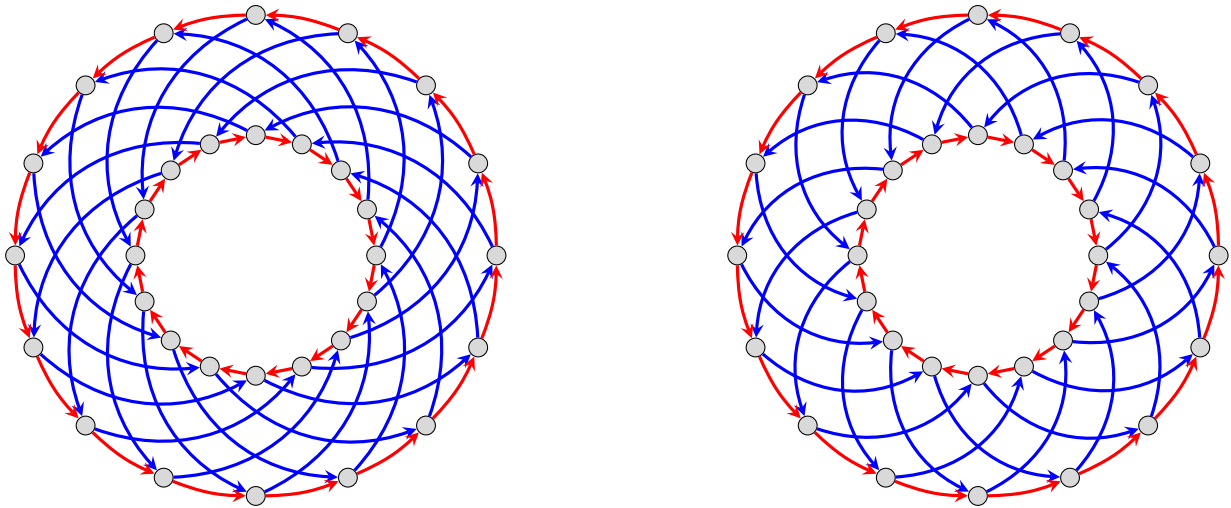

For another family of examples along the lines of the graphs in Figure~\ref{fig:flower-16}, consider the two graphs shown in Figure~\ref{fig:flower-32}. If both of these are Cayley graphs, then they must be either $D_{16}$, $\SA_{16}$, $\SA_{16}$, or $Q_{32}$. We encourage the reader to try to figure out if they are Cayley graphs of one of these, or if not, then what group they define by a presentation. These groups can be distinguished by several means. The dihedral group $D_{16}$ has $17$ elements of order $2$: all $16$ reflections, and the $180^\circ$ rotation. It and the semidihedral group $\SA_{16}$ both have a size-$2$ center, whereas $Z(\SA_{16})$ has four elements. The group $Q_{32}$ has only one element of order $2$.

%% Are these groups of order $32$?

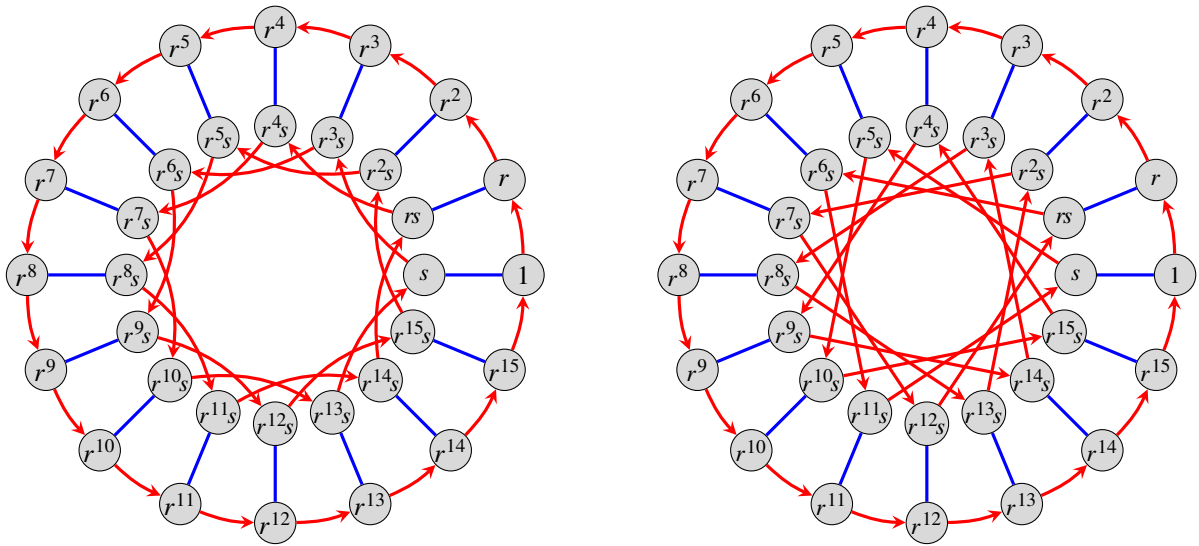
\begin{figure}[!ht] \centering
  \begin{tikzpicture}[scale=1.65,auto]
\tikzstyle{v} = [circle, draw, fill=lightgray,inner sep=0pt, minimum size=5.5mm]
\tikzstyle{r-out} = [draw, very thick, eRed,-stealth,bend right=10]
\tikzstyle{r-in} = [draw, very thick, eRed,-stealth,bend left=18]
    \tikzstyle{every node}=[font=\footnotesize]
    \begin{scope}[shift={(0,0)}]%,shorten >= -2pt, shorten <= -2pt]
      \node (s) at (0:1.2) [v] {$s$};
      \node (rs) at (22.5:1.2) [v] {$r\!s$};
      \node (r2s) at (45:1.2) [v] {$r^2\!s$};
      \node (r3s) at (67.5:1.2) [v] {$r^3\!s$};
      \node (r4s) at (90:1.2) [v] {$r^4\!s$};
      \node (r5s) at (112.5:1.2) [v] {$r^5\!s$};
      \node (r6s) at (135:1.2) [v] {$r^6\!s$};
      \node (r7s) at (157.5:1.2) [v] {$r^7\!s$};
      \node (r8s) at (180:1.2) [v] {$r^8\!s$};
      \node (r9s) at (-157.5:1.2) [v] {$r^9\!s$};
      \node (r10s) at (-135:1.2) [v] {$r^{10}\!s$};
      \node (r11s) at (-112.5:1.2) [v] {$r^{11}\!s$};
      \node (r12s) at (-90:1.2) [v] {$r^{12}\!s$};
      \node (r13s) at (-67.5:1.2) [v] {$r^{13}\!s$};
      \node (r14s) at (-45:1.2) [v] {$r^{14}\!s$};
      \node (r15s) at (-22.5:1.2) [v] {$r^{15}\!s$};
      \tikzstyle{every node}=[font=\small]
      \node (1) at (0:2) [v] {$1$};
      \node (r) at (22.5:2) [v] {$r$};
      \node (r2) at (45:2) [v] {$r^2$};
      \node (r3) at (67.5:2) [v] {$r^3$};
      \node (r4) at (90:2) [v] {$r^4$};
      \node (r5) at (112.5:2) [v] {$r^5$};
      \node (r6) at (135:2) [v] {$r^6$};
      \node (r7) at (157.5:2) [v] {$r^7$};
      \node (r8) at (180:2) [v] {$r^8$};
      \node (r9) at (-157.5:2) [v] {$r^9$};
      \node (r10) at (-135:2) [v] {$r^{10}$};
      \node (r11) at (-112.5:2) [v] {$r^{11}$};
      \node (r12) at (-90:2) [v] {$r^{12}$};
      \node (r13) at (-67.5:2) [v] {$r^{13}$};
      \node (r14) at (-45:2) [v] {$r^{14}$};
      \node (r15) at (-22.5:2) [v] {$r^{15}$};
      \draw [r-out] (1) to (r); \draw [r-out] (r) to (r2); \draw [r-out] (r2) to (r3);
      \draw [r-out] (r3) to (r4); \draw [r-out] (r4) to (r5); \draw [r-out] (r5) to (r6);
      \draw [r-out] (r6) to (r7); \draw [r-out] (r7) to (r8); \draw [r-out] (r8) to (r9);
      \draw [r-out] (r9) to (r10); \draw [r-out] (r10) to (r11);
      \draw [r-out] (r11) to (r12); \draw [r-out] (r12) to (r13);
      \draw [r-out] (r13) to (r14); \draw [r-out] (r14) to (r15);
      \draw [r-out] (r15) to (1);
      \draw [r-in] (s) to (r3s); \draw [r-in] (r3s) to (r6s);
      \draw [r-in] (r6s) to (r9s); \draw [r-in] (r9s) to (r12s);
      \draw [r-in] (r12s) to (r15s); \draw [r-in] (r15s) to (r2s);
      \draw [r-in] (r2s) to (r5s); \draw [r-in] (r5s) to (r8s);
      \draw [r-in] (r8s) to (r11s); \draw [r-in] (r11s) to (r14s);
      \draw [r-in] (r14s) to (rs); \draw [r-in] (rs) to (r4s);
      \draw [r-in] (r4s) to (r7s); \draw [r-in] (r7s) to (r10s);
      \draw [r-in] (r10s) to (r13s); \draw [r-in] (r13s) to (s); 
      \draw [bb] (1) to (s); \draw [bb] (r) to (rs);
      \draw [bb] (r2) to (r2s); \draw [bb] (r3) to (r3s);
      \draw [bb] (r4) to (r4s); \draw [bb] (r5) to (r5s);
      \draw [bb] (r6) to (r6s); \draw [bb] (r7) to (r7s);
      \draw [bb] (r8) to (r8s); \draw [bb] (r9) to (r9s);
      \draw [bb] (r10) to (r10s); \draw [bb] (r11) to (r11s);
      \draw [bb] (r12) to (r12s); \draw [bb] (r13) to (r13s);
      \draw [bb] (r14) to (r14s); \draw [bb] (r15) to (r15s);
    \end{scope}
    \begin{scope}[shift={(5.25,0)}]%,shorten >= -2pt, shorten <= -2pt]
      \node (s) at (0:1.2) [v] {$s$};
      \node (rs) at (22.5:1.2) [v] {$r\!s$};
      \node (r2s) at (45:1.2) [v] {$r^2\!s$};
      \node (r3s) at (67.5:1.2) [v] {$r^3\!s$};
      \node (r4s) at (90:1.2) [v] {$r^4\!s$};
      \node (r5s) at (112.5:1.2) [v] {$r^5\!s$};
      \node (r6s) at (135:1.2) [v] {$r^6\!s$};
      \node (r7s) at (157.5:1.2) [v] {$r^7\!s$};
      \node (r8s) at (180:1.2) [v] {$r^8\!s$};
      \node (r9s) at (-157.5:1.2) [v] {$r^9\!s$};
      \node (r10s) at (-135:1.2) [v] {$r^{10}\!s$};
      \node (r11s) at (-112.5:1.2) [v] {$r^{11}\!s$};
      \node (r12s) at (-90:1.2) [v] {$r^{12}\!s$};
      \node (r13s) at (-67.5:1.2) [v] {$r^{13}\!s$};
      \node (r14s) at (-45:1.2) [v] {$r^{14}\!s$};
      \node (r15s) at (-22.5:1.2) [v] {$r^{15}\!s$};
      \tikzstyle{every node}=[font=\footnotesize]
      \node (1) at (0:2) [v] {$1$};
      \node (r) at (22.5:2) [v] {$r$};
      \node (r2) at (45:2) [v] {$r^2$};
      \node (r3) at (67.5:2) [v] {$r^3$};
      \node (r4) at (90:2) [v] {$r^4$};
      \node (r5) at (112.5:2) [v] {$r^5$};
      \node (r6) at (135:2) [v] {$r^6$};
      \node (r7) at (157.5:2) [v] {$r^7$};
      \node (r8) at (180:2) [v] {$r^8$};
      \node (r9) at (-157.5:2) [v] {$r^9$};
      \node (r10) at (-135:2) [v] {$r^{10}$};
      \node (r11) at (-112.5:2) [v] {$r^{11}$};
      \node (r12) at (-90:2) [v] {$r^{12}$};
      \node (r13) at (-67.5:2) [v] {$r^{13}$};
      \node (r14) at (-45:2) [v] {$r^{14}$};
      \node (r15) at (-22.5:2) [v] {$r^{15}$};
      \draw [r-out] (1) to (r); \draw [r-out] (r) to (r2); \draw [r-out] (r2) to (r3);
      \draw [r-out] (r3) to (r4); \draw [r-out] (r4) to (r5); \draw [r-out] (r5) to (r6);
      \draw [r-out] (r6) to (r7); \draw [r-out] (r7) to (r8); \draw [r-out] (r8) to (r9);
      \draw [r-out] (r9) to (r10); \draw [r-out] (r10) to (r11);
      \draw [r-out] (r11) to (r12); \draw [r-out] (r12) to (r13);
      \draw [r-out] (r13) to (r14); \draw [r-out] (r14) to (r15);
      \draw [r-out] (r15) to (1);
      \draw [r] (s) to (r5s); \draw [r] (r5s) to (r10s);
      \draw [r] (r10s) to (r15s); \draw [r] (r15s) to (r4s);
      \draw [r] (r4s) to (r9s); \draw [r] (r9s) to (r14s);
      \draw [r] (r14s) to (r3s); \draw [r] (r3s) to (r8s);
      \draw [r] (r8s) to (r13s); \draw [r] (r13s) to (r2s);
      \draw [r] (r2s) to (r7s); \draw [r] (r7s) to (r12s);
      \draw [r] (r12s) to (rs); \draw [r] (rs) to (r6s);
      \draw [r] (r6s) to (r11s); \draw [r] (r11s) to (s); 
      \draw [bb] (1) to (s); \draw [bb] (r) to (rs);
      \draw [bb] (r2) to (r2s); \draw [bb] (r3) to (r3s);
      \draw [bb] (r4) to (r4s); \draw [bb] (r5) to (r5s);
      \draw [bb] (r6) to (r6s); \draw [bb] (r7) to (r7s);
      \draw [bb] (r8) to (r8s); \draw [bb] (r9) to (r9s);
      \draw [bb] (r10) to (r10s); \draw [bb] (r11) to (r11s);
      \draw [bb] (r12) to (r12s); \draw [bb] (r13) to (r13s);
      \draw [bb] (r14) to (r14s); \draw [bb] (r15) to (r15s);
    \end{scope}
  \end{tikzpicture}
  \caption{If these are Cayley graphs, then they must be either $D_{16}$, $\SD_{16}$, $\SA_{16}$, or $Q_{32}$.}\label{fig:32-nodes}
\end{figure}

Finally, we will conclude this section with one last pair of diagrams, shown in Figure~\ref{fig:32-nodes}. Recall that there are four semidirect products of $C_{16}$ with $C_2$, and each one is defined by the relation $srs=r^k$:  $k=1$ ($C_{16}\times C_2)$, $k=7$ ($\SD_{16}$), $k=9$ ($\SA_{16}$), and $k=15$ ($D_{16})$. The graph on the left satisfies the relation $srs=r^3$, and the one on the right, $srs=r^5$. Are these Cayley graphs for one of these four groups, or do they define presentations of smaller groups, like the Petersen graph example in Figure~\ref{fig:Petersen}?

%% 32 nodes, not a group?

\vspace{-5mm}

%%%%%%%%%%%%%
\section*{Summary and Conclusions}

There were four distinct goals of this paper. First and foremost, it is the aim of the author to dispel the notion that group theory is a dry and non-visual subject, as many students come away from a class thinking. Unfortunately, many algebra books do not promote this idea, though there are exceptions~\cite{carter2021visual}. Second, readers who teach abstract algebra might be inclined to incorporate some of these ideas, examples, and puzzles into their classes. Several of the groups in this article are simply inaccessible, or at least unmotivated, without the use of Cayley graphs. As a result, many students never see any finite groups beyond the ``usual suspects'' of the abelian, dihedral, symmetric, alternating, and quaternion groups, and these are limited in the types of behaviors that they exhibit. A third goal of this paper is to encourage interested readers to create new Cayley graph puzzles like the ones presented here. There are many ways to put together motifs to get graphs that ``look like'' Cayley graphs, but aren't. Finally, the last goal is to attract the interest of visual artists who enjoy incorporating mathematical patterns and symmetries into their work. There is an abundance of potential materials within the realm of group theory, much of which simply does not exist in the literature.

%%%%%%%%%%%%%%%%%%%%%%%%%%%%%%%%%%%%%%%
% References %
    
{\setlength{\baselineskip}{13pt} % tighten line spacing for bibliography
\raggedright				% no right justification for References
\bibliographystyle{bridges}

} % end setlength, raggedright
   
\end{document}